\documentclass{gtart}

\def\ifplaintex{\expandafter\ifx\csname documentclass\endcsname\relax}

\ifplaintex 
\else
\fi

\expandafter\ifx\csname beginpicture\endcsname\relax
\expandafter\ifx\csname documentclass\endcsname\relax
\input pictex \else
\input prepictex \input pictex \input postpictex \fi\fi

\def\gt{{\mathsurround=0pt\it $\cal G\mskip-2mu$eometry \&\ 
$\cal T\!\!$opology}}        

\def\gtp{{\mathsurround=0pt\it $\cal G\mskip-2mu$eometry \&\ 
$\cal T\!\!$opology $\cal P\!$ublications}}  

\def\lognumber#1{\def\thelognumber{#1}}
\def\volumenumber#1{\def\thevolumenumber{#1}}
\def\papernumber#1{\def\thepapernumber{#1}}
\def\volumeyear#1{\def\thevolumeyear{#1}}

\def\pagenumbers#1#2{\def\startpage{#1}\def\finishpage{#2}}
\def\published#1{\def\publishdate{#1}}
\def\proposed#1{\def\theproposer{#1}}
\def\seconded#1{\def\theseconders{#1}}
\def\received#1{\def\receiveddate{#1}}

\def\accepted#1{\def\accepteddate{#1}}

\def\asciiaddress#1{\def\theasciiaddress{#1}}
\def\asciiemail#1{\def\theasciiemail{#1}}
\long\def\asciiabstract#1{\long\def\theasciiabstract{#1}}
\def\asciikeywords#1{\def\theasciikeywords{#1}}

\let\\\par\let\thelognumber\relax
\let\thevolumenumber\relax\let\thepapernumber\relax
\let\thevolumeyear\relax\let\thesamplenumber\relax\let\startpage\relax
\let\finishpage\relax\let\publishdate\relax\let\receiveddate\relax
\let\reviseddate\relax\let\accepteddate\relax\let\theasciititle\relax
\let\theasciiauthors\relax\let\theasciiaddress\relax
\let\theasciiabstract\relax\let\theasciikeywords\relax
\let\theasciiemail\relax\let\theshortauthors\relax\let\theshorttitle\relax

\long\def\maketitlep{   

\count0=\startpage

\gt\hfill      
\beginpicture
\setcoordinatesystem units <0.33truein, 0.33truein> point at 2.2 0.9
\setplotsymbol ({$\cal G$})
\plotsymbolspacing=9truept
\circulararc 315 degrees from 0 1 center at 0 0
\setplotsymbol ({$\cal T$})
\circulararc 315 degrees from 1 -1 center at 1 0
\endpicture
\break
{\small\ifx\thesamplenumber\relax 
Volume \else Sample
\fi\thevolumenumber\ (\thevolumeyear)
\startpage--\finishpage\nl
Published: \publishdate}
\vglue 0.5truein plus 0.4fil minus 0.1truein

{\parskip=0pt\leftskip 0pt plus 1fil\def\\{\par\smallskip}{\ifplaintex\large
\else\Large\fi\bf\thetitle}\par\medskip}   

\vglue 0pt plus 0.1fil 

{\parskip=0pt\leftskip 0pt plus 1fil\def\\{\par}{\sc\theauthors}
\par\medskip}

\vglue 0pt plus 0.1fil 

{\small\parskip=0pt\let\newline\\
{\leftskip 0pt plus 1fil\def\\{\par}{\sl\theaddress}\par}
\expandafter\ifx\theemail\relax    
\relax\else\vglue 5pt plus 0.02fil minus 2pt\def\\{\stdspace{\rm 
and}\stdspace} 
\cl{Email:\stdspace\tt\theemail}\fi
\ifx\theurl\relax                  
\relax\else\vglue 5pt plus 0.02fil minus 2pt\def\\{\stdspace{\rm 
and}\stdspace}
\cl{URL:\stdspace\tt\theurl}\fi\par}

\vglue 7pt plus 0.3fil minus 3pt

{\bf Abstract}
\vglue 5pt plus 0.1fil minus 2pt

\theabstract

\vglue 7pt plus 0.3fil minus 3pt

{\bf AMS Classification numbers}\quad Primary:\quad \theprimaryclass

Secondary:\quad \thesecondaryclass

\vglue 5pt plus 0.3fil minus 2pt

{\bf Keywords}\quad \thekeywords

\vglue 10pt plus 0.5fil minus 5pt

{\small  Proposed: \theproposer\hfill Received: \receiveddate\nl
Seconded: \theseconders\hfill 
\ifx\reviseddate\relax                         
Accepted: \accepteddate                        
\else
Revised: \reviseddate                          
\fi}
\eject
}       

\let\maketitlepage\maketitlep
\let\maketitle\maketitlepage

\font\phead=cmsl9 scaled 950
\font=cmsl9 scaled 1050
\font\pnum=cmbx10 scaled 913
\font\lnum=cmbx10 
\font\pfoot=cmsl9 scaled 950
\font=cmsl9 scaled 1050
\ifplaintex
\headline{\vbox to 0pt{\vskip -4.5mm\line{\small\phead\ifnum
\count0=\startpage ISSN 1364-0380 (on line)
1465-3060 (printed) \hfill {\pnum\folio}\else\ifodd\count0\def\\{ }%
\ifx\theshorttitle\relax\thetitle\else\theshorttitle\fi\hfill{\pnum\folio}
\else\def\\{ and }{\pnum\folio}\hfill\ifx\theshortauthors\relax\theauthors
\else\theshortauthors\fi\fi\fi}\vss}}
\footline{\vbox to 0pt{\vglue 0mm\line{\small\pfoot\ifnum\count0=\startpage
\copyright\ \gtp\hfill\else
\gt, Volume \thevolumenumber\ (\thevolumeyear)\hfill\fi}\vss
}}
\else
\makeatletter
\def\@oddhead{{\small\lhead\ifnum\count0=\startpage ISSN 1364-0380 (on line)
1465-3060 (printed) \hfill {\lnum\number\count0}\else\ifodd\count0
\def\\{ }\ifx\theshorttitle\relax \thetitle \else\theshorttitle\fi\hfill
{\lnum\number\count0}\else\def\\{ and }{\lnum\number\count0}
\hfill\ifx\theshortauthors\relax 
\theauthors\else\theshortauthors\fi\fi\fi}}\def\@evenhead{\@oddhead}
\def\@oddfoot{\small\lfoot\ifnum\count0=\startpage\copyright\ \gtp\hfill\else
\gt, Volume \thevolumenumber\ (\thevolumeyear)\hfill\fi}
\def\@evenfoot{\@oddfoot}
\makeatother
\fi

\newwrite\gtoutfile
\long\gdef\makeheadfile{  
{\def\\{, }\def\s{ }
\immediate\openout\gtoutfile head.xxx
\immediate\write\gtoutfile{To: math@arxiv.org}
\immediate\write\gtoutfile{Subject: put or rep NNNNN:pppp}
\immediate\write\gtoutfile{--text follows this line--}
\immediate\write\gtoutfile{Proxy-for: \ifx\theasciiauthors\relax
\theauthors\else\theasciiauthors\fi\s<\ifx\theasciiemail\relax\theemail\else\theasciiemail\fi>}
\immediate\write\gtoutfile{\noexpand\\}
\immediate\write\gtoutfile{Authors: \ifx\theasciiauthors\relax
\theauthors\else\theasciiauthors\fi}
\immediate\write\gtoutfile{Title: \ifx\theasciititle\relax
\thetitle\else\theasciititle\fi}
\immediate\write\gtoutfile{Subj-class: GT or SG or MG etc}
\immediate\write\gtoutfile{MSC-class: \theprimaryclass\ifx\thesecondaryclass\relax\else, \thesecondaryclass\fi}
\immediate\write\gtoutfile{Journal-ref: Geom. Topol. \thevolumenumber
(\thevolumeyear) \startpage-\finishpage}
\immediate\write\gtoutfile{Comments: Published by Geometry and Topology at}
\immediate\write\gtoutfile{\s\s http://www.maths.warwick.ac.uk/gt/GTVol\thevolumenumber/paper\thepapernumber.abs.html}
\immediate\write\gtoutfile{\noexpand\\}
\immediate\write\gtoutfile{}
\ifx\theasciiabstract\relax
\immediate\write\gtoutfile{\theabstract}\else
\immediate\write\gtoutfile{\theasciiabstract}\fi
\immediate\write\gtoutfile{}
\immediate\write\gtoutfile{\noexpand\\}
\immediate\write\gtoutfile{}
\immediate\closeout\gtoutfile}}  

\def\maketitlepage{\maketitlep\makeheadfile}
\let\maketitle\maketitlepage

\def\ifplaintex{\expandafter\ifx\csname documentclass\endcsname\relax}

\ifplaintex 
\else
\fi

\expandafter\ifx\csname beginpicture\endcsname\relax
\expandafter\ifx\csname documentclass\endcsname\relax
\input pictex \else
\input prepictex \input pictex \input postpictex \fi\fi

\def\gt{{\mathsurround=0pt\it $\cal G\mskip-2mu$eometry \&\ 
$\cal T\!\!$opology}}        

\def\gtp{{\mathsurround=0pt\it $\cal G\mskip-2mu$eometry \&\ 
$\cal T\!\!$opology $\cal P\!$ublications}}  

\def\lognumber#1{\def\thelognumber{#1}}
\def\volumenumber#1{\def\thevolumenumber{#1}}
\def\papernumber#1{\def\thepapernumber{#1}}
\def\volumeyear#1{\def\thevolumeyear{#1}}

\def\pagenumbers#1#2{\def\startpage{#1}\def\finishpage{#2}}
\def\published#1{\def\publishdate{#1}}
\def\proposed#1{\def\theproposer{#1}}
\def\seconded#1{\def\theseconders{#1}}
\def\received#1{\def\receiveddate{#1}}

\def\accepted#1{\def\accepteddate{#1}}

\def\asciiaddress#1{\def\theasciiaddress{#1}}
\def\asciiemail#1{\def\theasciiemail{#1}}
\long\def\asciiabstract#1{\long\def\theasciiabstract{#1}}
\def\asciikeywords#1{\def\theasciikeywords{#1}}

\let\\\par\let\thelognumber\relax
\let\thevolumenumber\relax\let\thepapernumber\relax
\let\thevolumeyear\relax\let\thesamplenumber\relax\let\startpage\relax
\let\finishpage\relax\let\publishdate\relax\let\receiveddate\relax
\let\reviseddate\relax\let\accepteddate\relax\let\theasciititle\relax
\let\theasciiauthors\relax\let\theasciiaddress\relax
\let\theasciiabstract\relax\let\theasciikeywords\relax
\let\theasciiemail\relax\let\theshortauthors\relax\let\theshorttitle\relax

\long\def\maketitlep{   

\count0=\startpage

\gt\hfill      
\beginpicture
\setcoordinatesystem units <0.33truein, 0.33truein> point at 2.2 0.9
\setplotsymbol ({$\cal G$})
\plotsymbolspacing=9truept
\circulararc 315 degrees from 0 1 center at 0 0
\setplotsymbol ({$\cal T$})
\circulararc 315 degrees from 1 -1 center at 1 0
\endpicture
\break
{\small\ifx\thesamplenumber\relax 
Volume \else Sample
\fi\thevolumenumber\ (\thevolumeyear)
\startpage--\finishpage\nl
Published: \publishdate}
\vglue 0.5truein plus 0.4fil minus 0.1truein

{\parskip=0pt\leftskip 0pt plus 1fil\def\\{\par\smallskip}{\ifplaintex\large
\else\Large\fi\bf\thetitle}\par\medskip}   

\vglue 0pt plus 0.1fil 

{\parskip=0pt\leftskip 0pt plus 1fil\def\\{\par}{\sc\theauthors}
\par\medskip}

\vglue 0pt plus 0.1fil 

{\small\parskip=0pt\let\newline\\
{\leftskip 0pt plus 1fil\def\\{\par}{\sl\theaddress}\par}
\expandafter\ifx\theemail\relax    
\relax\else\vglue 5pt plus 0.02fil minus 2pt\def\\{\stdspace{\rm 
and}\stdspace} 
\cl{Email:\stdspace\tt\theemail}\fi
\ifx\theurl\relax                  
\relax\else\vglue 5pt plus 0.02fil minus 2pt\def\\{\stdspace{\rm 
and}\stdspace}
\cl{URL:\stdspace\tt\theurl}\fi\par}

\vglue 7pt plus 0.3fil minus 3pt

{\bf Abstract}
\vglue 5pt plus 0.1fil minus 2pt

\theabstract

\vglue 7pt plus 0.3fil minus 3pt

{\bf AMS Classification numbers}\quad Primary:\quad \theprimaryclass

Secondary:\quad \thesecondaryclass

\vglue 5pt plus 0.3fil minus 2pt

{\bf Keywords}\quad \thekeywords

\vglue 10pt plus 0.5fil minus 5pt

{\small  Proposed: \theproposer\hfill Received: \receiveddate\nl
Seconded: \theseconders\hfill 
\ifx\reviseddate\relax                         
Accepted: \accepteddate                        
\else
Revised: \reviseddate                          
\fi}
\eject
}       

\let\maketitlepage\maketitlep
\let\maketitle\maketitlepage

\font\phead=cmsl9 scaled 950
\font=cmsl9 scaled 1050
\font\pnum=cmbx10 scaled 913
\font\lnum=cmbx10 
\font\pfoot=cmsl9 scaled 950
\font=cmsl9 scaled 1050
\ifplaintex
\headline{\vbox to 0pt{\vskip -4.5mm\line{\small\phead\ifnum
\count0=\startpage ISSN 1364-0380 (on line)
1465-3060 (printed) \hfill {\pnum\folio}\else\ifodd\count0\def\\{ }%
\ifx\theshorttitle\relax\thetitle\else\theshorttitle\fi\hfill{\pnum\folio}
\else\def\\{ and }{\pnum\folio}\hfill\ifx\theshortauthors\relax\theauthors
\else\theshortauthors\fi\fi\fi}\vss}}
\footline{\vbox to 0pt{\vglue 0mm\line{\small\pfoot\ifnum\count0=\startpage
\copyright\ \gtp\hfill\else
\gt, Volume \thevolumenumber\ (\thevolumeyear)\hfill\fi}\vss
}}
\else
\makeatletter
\def\@oddhead{{\small\lhead\ifnum\count0=\startpage ISSN 1364-0380 (on line)
1465-3060 (printed) \hfill {\lnum\number\count0}\else\ifodd\count0
\def\\{ }\ifx\theshorttitle\relax \thetitle \else\theshorttitle\fi\hfill
{\lnum\number\count0}\else\def\\{ and }{\lnum\number\count0}
\hfill\ifx\theshortauthors\relax 
\theauthors\else\theshortauthors\fi\fi\fi}}\def\@evenhead{\@oddhead}
\def\@oddfoot{\small\lfoot\ifnum\count0=\startpage\copyright\ \gtp\hfill\else
\gt, Volume \thevolumenumber\ (\thevolumeyear)\hfill\fi}
\def\@evenfoot{\@oddfoot}
\makeatother
\fi

\newwrite\gtoutfile
\long\gdef\makeheadfile{  
{\def\\{, }\def\s{ }
\immediate\openout\gtoutfile head.xxx
\immediate\write\gtoutfile{To: math@arxiv.org}
\immediate\write\gtoutfile{Subject: put or rep NNNNN:pppp}
\immediate\write\gtoutfile{--text follows this line--}
\immediate\write\gtoutfile{Proxy-for: \ifx\theasciiauthors\relax
\theauthors\else\theasciiauthors\fi\s<\ifx\theasciiemail\relax\theemail\else\theasciiemail\fi>}
\immediate\write\gtoutfile{\noexpand\\}
\immediate\write\gtoutfile{Authors: \ifx\theasciiauthors\relax
\theauthors\else\theasciiauthors\fi}
\immediate\write\gtoutfile{Title: \ifx\theasciititle\relax
\thetitle\else\theasciititle\fi}
\immediate\write\gtoutfile{Subj-class: GT or SG or MG etc}
\immediate\write\gtoutfile{MSC-class: \theprimaryclass\ifx\thesecondaryclass\relax\else, \thesecondaryclass\fi}
\immediate\write\gtoutfile{Journal-ref: Geom. Topol. \thevolumenumber
(\thevolumeyear) \startpage-\finishpage}
\immediate\write\gtoutfile{Comments: Published by Geometry and Topology at}
\immediate\write\gtoutfile{\s\s http://www.maths.warwick.ac.uk/gt/GTVol\thevolumenumber/paper\thepapernumber.abs.html}
\immediate\write\gtoutfile{\noexpand\\}
\immediate\write\gtoutfile{}
\ifx\theasciiabstract\relax
\immediate\write\gtoutfile{\theabstract}\else
\immediate\write\gtoutfile{\theasciiabstract}\fi
\immediate\write\gtoutfile{}
\immediate\write\gtoutfile{\noexpand\\}
\immediate\write\gtoutfile{}
\immediate\closeout\gtoutfile}}  

\def\maketitlepage{\maketitlep\makeheadfile}
\let\maketitle\maketitlepage

\def\ifplaintex{\expandafter\ifx\csname documentclass\endcsname\relax}

\ifplaintex 
\else
\fi

\expandafter\ifx\csname beginpicture\endcsname\relax
\expandafter\ifx\csname documentclass\endcsname\relax
\input pictex \else
\input prepictex \input pictex \input postpictex \fi\fi

\def\gt{{\mathsurround=0pt\it $\cal G\mskip-2mu$eometry \&\ 
$\cal T\!\!$opology}}        

\def\gtp{{\mathsurround=0pt\it $\cal G\mskip-2mu$eometry \&\ 
$\cal T\!\!$opology $\cal P\!$ublications}}  

\def\lognumber#1{\def\thelognumber{#1}}
\def\volumenumber#1{\def\thevolumenumber{#1}}
\def\papernumber#1{\def\thepapernumber{#1}}
\def\volumeyear#1{\def\thevolumeyear{#1}}

\def\pagenumbers#1#2{\def\startpage{#1}\def\finishpage{#2}}
\def\published#1{\def\publishdate{#1}}
\def\proposed#1{\def\theproposer{#1}}
\def\seconded#1{\def\theseconders{#1}}
\def\received#1{\def\receiveddate{#1}}

\def\accepted#1{\def\accepteddate{#1}}

\def\asciiaddress#1{\def\theasciiaddress{#1}}
\def\asciiemail#1{\def\theasciiemail{#1}}
\long\def\asciiabstract#1{\long\def\theasciiabstract{#1}}
\def\asciikeywords#1{\def\theasciikeywords{#1}}

\let\\\par\let\thelognumber\relax
\let\thevolumenumber\relax\let\thepapernumber\relax
\let\thevolumeyear\relax\let\thesamplenumber\relax\let\startpage\relax
\let\finishpage\relax\let\publishdate\relax\let\receiveddate\relax
\let\reviseddate\relax\let\accepteddate\relax\let\theasciititle\relax
\let\theasciiauthors\relax\let\theasciiaddress\relax
\let\theasciiabstract\relax\let\theasciikeywords\relax
\let\theasciiemail\relax\let\theshortauthors\relax\let\theshorttitle\relax

\long\def\maketitlep{   

\count0=\startpage

\gt\hfill      
\beginpicture
\setcoordinatesystem units <0.33truein, 0.33truein> point at 2.2 0.9
\setplotsymbol ({$\cal G$})
\plotsymbolspacing=9truept
\circulararc 315 degrees from 0 1 center at 0 0
\setplotsymbol ({$\cal T$})
\circulararc 315 degrees from 1 -1 center at 1 0
\endpicture
\break
{\small\ifx\thesamplenumber\relax 
Volume \else Sample
\fi\thevolumenumber\ (\thevolumeyear)
\startpage--\finishpage\nl
Published: \publishdate}
\vglue 0.5truein plus 0.4fil minus 0.1truein

{\parskip=0pt\leftskip 0pt plus 1fil\def\\{\par\smallskip}{\ifplaintex\large
\else\Large\fi\bf\thetitle}\par\medskip}   

\vglue 0pt plus 0.1fil 

{\parskip=0pt\leftskip 0pt plus 1fil\def\\{\par}{\sc\theauthors}
\par\medskip}

\vglue 0pt plus 0.1fil 

{\small\parskip=0pt\let\newline\\
{\leftskip 0pt plus 1fil\def\\{\par}{\sl\theaddress}\par}
\expandafter\ifx\theemail\relax    
\relax\else\vglue 5pt plus 0.02fil minus 2pt\def\\{\stdspace{\rm 
and}\stdspace} 
\cl{Email:\stdspace\tt\theemail}\fi
\ifx\theurl\relax                  
\relax\else\vglue 5pt plus 0.02fil minus 2pt\def\\{\stdspace{\rm 
and}\stdspace}
\cl{URL:\stdspace\tt\theurl}\fi\par}

\vglue 7pt plus 0.3fil minus 3pt

{\bf Abstract}
\vglue 5pt plus 0.1fil minus 2pt

\theabstract

\vglue 7pt plus 0.3fil minus 3pt

{\bf AMS Classification numbers}\quad Primary:\quad \theprimaryclass

Secondary:\quad \thesecondaryclass

\vglue 5pt plus 0.3fil minus 2pt

{\bf Keywords}\quad \thekeywords

\vglue 10pt plus 0.5fil minus 5pt

{\small  Proposed: \theproposer\hfill Received: \receiveddate\nl
Seconded: \theseconders\hfill 
\ifx\reviseddate\relax                         
Accepted: \accepteddate                        
\else
Revised: \reviseddate                          
\fi}
\eject
}       

\let\maketitlepage\maketitlep
\let\maketitle\maketitlepage

\font\phead=cmsl9 scaled 950
\font=cmsl9 scaled 1050
\font\pnum=cmbx10 scaled 913
\font\lnum=cmbx10 
\font\pfoot=cmsl9 scaled 950
\font=cmsl9 scaled 1050
\ifplaintex
\headline{\vbox to 0pt{\vskip -4.5mm\line{\small\phead\ifnum
\count0=\startpage ISSN 1364-0380 (on line)
1465-3060 (printed) \hfill {\pnum\folio}\else\ifodd\count0\def\\{ }%
\ifx\theshorttitle\relax\thetitle\else\theshorttitle\fi\hfill{\pnum\folio}
\else\def\\{ and }{\pnum\folio}\hfill\ifx\theshortauthors\relax\theauthors
\else\theshortauthors\fi\fi\fi}\vss}}
\footline{\vbox to 0pt{\vglue 0mm\line{\small\pfoot\ifnum\count0=\startpage
\copyright\ \gtp\hfill\else
\gt, Volume \thevolumenumber\ (\thevolumeyear)\hfill\fi}\vss
}}
\else
\makeatletter
\def\@oddhead{{\small\lhead\ifnum\count0=\startpage ISSN 1364-0380 (on line)
1465-3060 (printed) \hfill {\lnum\number\count0}\else\ifodd\count0
\def\\{ }\ifx\theshorttitle\relax \thetitle \else\theshorttitle\fi\hfill
{\lnum\number\count0}\else\def\\{ and }{\lnum\number\count0}
\hfill\ifx\theshortauthors\relax 
\theauthors\else\theshortauthors\fi\fi\fi}}\def\@evenhead{\@oddhead}
\def\@oddfoot{\small\lfoot\ifnum\count0=\startpage\copyright\ \gtp\hfill\else
\gt, Volume \thevolumenumber\ (\thevolumeyear)\hfill\fi}
\def\@evenfoot{\@oddfoot}
\makeatother
\fi

\newwrite\gtoutfile
\long\gdef\makeheadfile{  
{\def\\{, }\def\s{ }
\immediate\openout\gtoutfile head.xxx
\immediate\write\gtoutfile{To: math@arxiv.org}
\immediate\write\gtoutfile{Subject: put or rep NNNNN:pppp}
\immediate\write\gtoutfile{--text follows this line--}
\immediate\write\gtoutfile{Proxy-for: \ifx\theasciiauthors\relax
\theauthors\else\theasciiauthors\fi\s<\ifx\theasciiemail\relax\theemail\else\theasciiemail\fi>}
\immediate\write\gtoutfile{\noexpand\\}
\immediate\write\gtoutfile{Authors: \ifx\theasciiauthors\relax
\theauthors\else\theasciiauthors\fi}
\immediate\write\gtoutfile{Title: \ifx\theasciititle\relax
\thetitle\else\theasciititle\fi}
\immediate\write\gtoutfile{Subj-class: GT or SG or MG etc}
\immediate\write\gtoutfile{MSC-class: \theprimaryclass\ifx\thesecondaryclass\relax\else, \thesecondaryclass\fi}
\immediate\write\gtoutfile{Journal-ref: Geom. Topol. \thevolumenumber
(\thevolumeyear) \startpage-\finishpage}
\immediate\write\gtoutfile{Comments: Published by Geometry and Topology at}
\immediate\write\gtoutfile{\s\s http://www.maths.warwick.ac.uk/gt/GTVol\thevolumenumber/paper\thepapernumber.abs.html}
\immediate\write\gtoutfile{\noexpand\\}
\immediate\write\gtoutfile{}
\ifx\theasciiabstract\relax
\immediate\write\gtoutfile{\theabstract}\else
\immediate\write\gtoutfile{\theasciiabstract}\fi
\immediate\write\gtoutfile{}
\immediate\write\gtoutfile{\noexpand\\}
\immediate\write\gtoutfile{}
\immediate\closeout\gtoutfile}}  

\def\maketitlepage{\maketitlep\makeheadfile}
\let\maketitle\maketitlepage

\def\ifplaintex{\expandafter\ifx\csname documentclass\endcsname\relax}

\ifplaintex 
\else
\fi

\expandafter\ifx\csname beginpicture\endcsname\relax
\expandafter\ifx\csname documentclass\endcsname\relax
\input pictex \else
\input prepictex \input pictex \input postpictex \fi\fi

\def\gt{{\mathsurround=0pt\it $\cal G\mskip-2mu$eometry \&\ 
$\cal T\!\!$opology}}        

\def\gtp{{\mathsurround=0pt\it $\cal G\mskip-2mu$eometry \&\ 
$\cal T\!\!$opology $\cal P\!$ublications}}  

\def\lognumber#1{\def\thelognumber{#1}}
\def\volumenumber#1{\def\thevolumenumber{#1}}
\def\papernumber#1{\def\thepapernumber{#1}}
\def\volumeyear#1{\def\thevolumeyear{#1}}

\def\pagenumbers#1#2{\def\startpage{#1}\def\finishpage{#2}}
\def\published#1{\def\publishdate{#1}}
\def\proposed#1{\def\theproposer{#1}}
\def\seconded#1{\def\theseconders{#1}}
\def\received#1{\def\receiveddate{#1}}

\def\accepted#1{\def\accepteddate{#1}}

\def\asciiaddress#1{\def\theasciiaddress{#1}}
\def\asciiemail#1{\def\theasciiemail{#1}}
\long\def\asciiabstract#1{\long\def\theasciiabstract{#1}}
\def\asciikeywords#1{\def\theasciikeywords{#1}}

\let\\\par\let\thelognumber\relax
\let\thevolumenumber\relax\let\thepapernumber\relax
\let\thevolumeyear\relax\let\thesamplenumber\relax\let\startpage\relax
\let\finishpage\relax\let\publishdate\relax\let\receiveddate\relax
\let\reviseddate\relax\let\accepteddate\relax\let\theasciititle\relax
\let\theasciiauthors\relax\let\theasciiaddress\relax
\let\theasciiabstract\relax\let\theasciikeywords\relax
\let\theasciiemail\relax\let\theshortauthors\relax\let\theshorttitle\relax

\long\def\maketitlep{   

\count0=\startpage

\gt\hfill      
\beginpicture
\setcoordinatesystem units <0.33truein, 0.33truein> point at 2.2 0.9
\setplotsymbol ({$\cal G$})
\plotsymbolspacing=9truept
\circulararc 315 degrees from 0 1 center at 0 0
\setplotsymbol ({$\cal T$})
\circulararc 315 degrees from 1 -1 center at 1 0
\endpicture
\break
{\small\ifx\thesamplenumber\relax 
Volume \else Sample
\fi\thevolumenumber\ (\thevolumeyear)
\startpage--\finishpage\nl
Published: \publishdate}
\vglue 0.5truein plus 0.4fil minus 0.1truein

{\parskip=0pt\leftskip 0pt plus 1fil\def\\{\par\smallskip}{\ifplaintex\large
\else\Large\fi\bf\thetitle}\par\medskip}   

\vglue 0pt plus 0.1fil 

{\parskip=0pt\leftskip 0pt plus 1fil\def\\{\par}{\sc\theauthors}
\par\medskip}

\vglue 0pt plus 0.1fil 

{\small\parskip=0pt\let\newline\\
{\leftskip 0pt plus 1fil\def\\{\par}{\sl\theaddress}\par}
\expandafter\ifx\theemail\relax    
\relax\else\vglue 5pt plus 0.02fil minus 2pt\def\\{\stdspace{\rm 
and}\stdspace} 
\cl{Email:\stdspace\tt\theemail}\fi
\ifx\theurl\relax                  
\relax\else\vglue 5pt plus 0.02fil minus 2pt\def\\{\stdspace{\rm 
and}\stdspace}
\cl{URL:\stdspace\tt\theurl}\fi\par}

\vglue 7pt plus 0.3fil minus 3pt

{\bf Abstract}
\vglue 5pt plus 0.1fil minus 2pt

\theabstract

\vglue 7pt plus 0.3fil minus 3pt

{\bf AMS Classification numbers}\quad Primary:\quad \theprimaryclass

Secondary:\quad \thesecondaryclass

\vglue 5pt plus 0.3fil minus 2pt

{\bf Keywords}\quad \thekeywords

\vglue 10pt plus 0.5fil minus 5pt

{\small  Proposed: \theproposer\hfill Received: \receiveddate\nl
Seconded: \theseconders\hfill 
\ifx\reviseddate\relax                         
Accepted: \accepteddate                        
\else
Revised: \reviseddate                          
\fi}
\eject
}       

\let\maketitlepage\maketitlep
\let\maketitle\maketitlepage

\font\phead=cmsl9 scaled 950
\font=cmsl9 scaled 1050
\font\pnum=cmbx10 scaled 913
\font\lnum=cmbx10 
\font\pfoot=cmsl9 scaled 950
\font=cmsl9 scaled 1050
\ifplaintex
\headline{\vbox to 0pt{\vskip -4.5mm\line{\small\phead\ifnum
\count0=\startpage ISSN 1364-0380 (on line)
1465-3060 (printed) \hfill {\pnum\folio}\else\ifodd\count0\def\\{ }%
\ifx\theshorttitle\relax\thetitle\else\theshorttitle\fi\hfill{\pnum\folio}
\else\def\\{ and }{\pnum\folio}\hfill\ifx\theshortauthors\relax\theauthors
\else\theshortauthors\fi\fi\fi}\vss}}
\footline{\vbox to 0pt{\vglue 0mm\line{\small\pfoot\ifnum\count0=\startpage
\copyright\ \gtp\hfill\else
\gt, Volume \thevolumenumber\ (\thevolumeyear)\hfill\fi}\vss
}}
\else
\makeatletter
\def\@oddhead{{\small\lhead\ifnum\count0=\startpage ISSN 1364-0380 (on line)
1465-3060 (printed) \hfill {\lnum\number\count0}\else\ifodd\count0
\def\\{ }\ifx\theshorttitle\relax \thetitle \else\theshorttitle\fi\hfill
{\lnum\number\count0}\else\def\\{ and }{\lnum\number\count0}
\hfill\ifx\theshortauthors\relax 
\theauthors\else\theshortauthors\fi\fi\fi}}\def\@evenhead{\@oddhead}
\def\@oddfoot{\small\lfoot\ifnum\count0=\startpage\copyright\ \gtp\hfill\else
\gt, Volume \thevolumenumber\ (\thevolumeyear)\hfill\fi}
\def\@evenfoot{\@oddfoot}
\makeatother
\fi

\newwrite\gtoutfile
\long\gdef\makeheadfile{  
{\def\\{, }\def\s{ }
\immediate\openout\gtoutfile head.xxx
\immediate\write\gtoutfile{To: math@arxiv.org}
\immediate\write\gtoutfile{Subject: put or rep NNNNN:pppp}
\immediate\write\gtoutfile{--text follows this line--}
\immediate\write\gtoutfile{Proxy-for: \ifx\theasciiauthors\relax
\theauthors\else\theasciiauthors\fi\s<\ifx\theasciiemail\relax\theemail\else\theasciiemail\fi>}
\immediate\write\gtoutfile{\noexpand\\}
\immediate\write\gtoutfile{Authors: \ifx\theasciiauthors\relax
\theauthors\else\theasciiauthors\fi}
\immediate\write\gtoutfile{Title: \ifx\theasciititle\relax
\thetitle\else\theasciititle\fi}
\immediate\write\gtoutfile{Subj-class: GT or SG or MG etc}
\immediate\write\gtoutfile{MSC-class: \theprimaryclass\ifx\thesecondaryclass\relax\else, \thesecondaryclass\fi}
\immediate\write\gtoutfile{Journal-ref: Geom. Topol. \thevolumenumber
(\thevolumeyear) \startpage-\finishpage}
\immediate\write\gtoutfile{Comments: Published by Geometry and Topology at}
\immediate\write\gtoutfile{\s\s http://www.maths.warwick.ac.uk/gt/GTVol\thevolumenumber/paper\thepapernumber.abs.html}
\immediate\write\gtoutfile{\noexpand\\}
\immediate\write\gtoutfile{}
\ifx\theasciiabstract\relax
\immediate\write\gtoutfile{\theabstract}\else
\immediate\write\gtoutfile{\theasciiabstract}\fi
\immediate\write\gtoutfile{}
\immediate\write\gtoutfile{\noexpand\\}
\immediate\write\gtoutfile{}
\immediate\closeout\gtoutfile}}  

\def\maketitlepage{\maketitlep\makeheadfile}
\let\maketitle\maketitlepage

\def\ifplaintex{\expandafter\ifx\csname documentclass\endcsname\relax}

\ifplaintex 
\else
\fi

\expandafter\ifx\csname beginpicture\endcsname\relax
\expandafter\ifx\csname documentclass\endcsname\relax
\input pictex \else
\input prepictex \input pictex \input postpictex \fi\fi

\def\gt{{\mathsurround=0pt\it $\cal G\mskip-2mu$eometry \&\ 
$\cal T\!\!$opology}}        

\def\gtp{{\mathsurround=0pt\it $\cal G\mskip-2mu$eometry \&\ 
$\cal T\!\!$opology $\cal P\!$ublications}}  

\def\lognumber#1{\def\thelognumber{#1}}
\def\volumenumber#1{\def\thevolumenumber{#1}}
\def\papernumber#1{\def\thepapernumber{#1}}
\def\volumeyear#1{\def\thevolumeyear{#1}}

\def\pagenumbers#1#2{\def\startpage{#1}\def\finishpage{#2}}
\def\published#1{\def\publishdate{#1}}
\def\proposed#1{\def\theproposer{#1}}
\def\seconded#1{\def\theseconders{#1}}
\def\received#1{\def\receiveddate{#1}}

\def\accepted#1{\def\accepteddate{#1}}

\def\asciiaddress#1{\def\theasciiaddress{#1}}
\def\asciiemail#1{\def\theasciiemail{#1}}
\long\def\asciiabstract#1{\long\def\theasciiabstract{#1}}
\def\asciikeywords#1{\def\theasciikeywords{#1}}

\let\\\par\let\thelognumber\relax
\let\thevolumenumber\relax\let\thepapernumber\relax
\let\thevolumeyear\relax\let\thesamplenumber\relax\let\startpage\relax
\let\finishpage\relax\let\publishdate\relax\let\receiveddate\relax
\let\reviseddate\relax\let\accepteddate\relax\let\theasciititle\relax
\let\theasciiauthors\relax\let\theasciiaddress\relax
\let\theasciiabstract\relax\let\theasciikeywords\relax
\let\theasciiemail\relax\let\theshortauthors\relax\let\theshorttitle\relax

\long\def\maketitlep{   

\count0=\startpage

\gt\hfill      
\beginpicture
\setcoordinatesystem units <0.33truein, 0.33truein> point at 2.2 0.9
\setplotsymbol ({$\cal G$})
\plotsymbolspacing=9truept
\circulararc 315 degrees from 0 1 center at 0 0
\setplotsymbol ({$\cal T$})
\circulararc 315 degrees from 1 -1 center at 1 0
\endpicture
\break
{\small\ifx\thesamplenumber\relax 
Volume \else Sample
\fi\thevolumenumber\ (\thevolumeyear)
\startpage--\finishpage\nl
Published: \publishdate}
\vglue 0.5truein plus 0.4fil minus 0.1truein

{\parskip=0pt\leftskip 0pt plus 1fil\def\\{\par\smallskip}{\ifplaintex\large
\else\Large\fi\bf\thetitle}\par\medskip}   

\vglue 0pt plus 0.1fil 

{\parskip=0pt\leftskip 0pt plus 1fil\def\\{\par}{\sc\theauthors}
\par\medskip}

\vglue 0pt plus 0.1fil 

{\small\parskip=0pt\let\newline\\
{\leftskip 0pt plus 1fil\def\\{\par}{\sl\theaddress}\par}
\expandafter\ifx\theemail\relax    
\relax\else\vglue 5pt plus 0.02fil minus 2pt\def\\{\stdspace{\rm 
and}\stdspace} 
\cl{Email:\stdspace\tt\theemail}\fi
\ifx\theurl\relax                  
\relax\else\vglue 5pt plus 0.02fil minus 2pt\def\\{\stdspace{\rm 
and}\stdspace}
\cl{URL:\stdspace\tt\theurl}\fi\par}

\vglue 7pt plus 0.3fil minus 3pt

{\bf Abstract}
\vglue 5pt plus 0.1fil minus 2pt

\theabstract

\vglue 7pt plus 0.3fil minus 3pt

{\bf AMS Classification numbers}\quad Primary:\quad \theprimaryclass

Secondary:\quad \thesecondaryclass

\vglue 5pt plus 0.3fil minus 2pt

{\bf Keywords}\quad \thekeywords

\vglue 10pt plus 0.5fil minus 5pt

{\small  Proposed: \theproposer\hfill Received: \receiveddate\nl
Seconded: \theseconders\hfill 
\ifx\reviseddate\relax                         
Accepted: \accepteddate                        
\else
Revised: \reviseddate                          
\fi}
\eject
}       

\let\maketitlepage\maketitlep
\let\maketitle\maketitlepage

\font\phead=cmsl9 scaled 950
\font=cmsl9 scaled 1050
\font\pnum=cmbx10 scaled 913
\font\lnum=cmbx10 
\font\pfoot=cmsl9 scaled 950
\font=cmsl9 scaled 1050
\ifplaintex
\headline{\vbox to 0pt{\vskip -4.5mm\line{\small\phead\ifnum
\count0=\startpage ISSN 1364-0380 (on line)
1465-3060 (printed) \hfill {\pnum\folio}\else\ifodd\count0\def\\{ }%
\ifx\theshorttitle\relax\thetitle\else\theshorttitle\fi\hfill{\pnum\folio}
\else\def\\{ and }{\pnum\folio}\hfill\ifx\theshortauthors\relax\theauthors
\else\theshortauthors\fi\fi\fi}\vss}}
\footline{\vbox to 0pt{\vglue 0mm\line{\small\pfoot\ifnum\count0=\startpage
\copyright\ \gtp\hfill\else
\gt, Volume \thevolumenumber\ (\thevolumeyear)\hfill\fi}\vss
}}
\else
\makeatletter
\def\@oddhead{{\small\lhead\ifnum\count0=\startpage ISSN 1364-0380 (on line)
1465-3060 (printed) \hfill {\lnum\number\count0}\else\ifodd\count0
\def\\{ }\ifx\theshorttitle\relax \thetitle \else\theshorttitle\fi\hfill
{\lnum\number\count0}\else\def\\{ and }{\lnum\number\count0}
\hfill\ifx\theshortauthors\relax 
\theauthors\else\theshortauthors\fi\fi\fi}}\def\@evenhead{\@oddhead}
\def\@oddfoot{\small\lfoot\ifnum\count0=\startpage\copyright\ \gtp\hfill\else
\gt, Volume \thevolumenumber\ (\thevolumeyear)\hfill\fi}
\def\@evenfoot{\@oddfoot}
\makeatother
\fi

\newwrite\gtoutfile
\long\gdef\makeheadfile{  
{\def\\{, }\def\s{ }
\immediate\openout\gtoutfile head.xxx
\immediate\write\gtoutfile{To: math@arxiv.org}
\immediate\write\gtoutfile{Subject: put or rep NNNNN:pppp}
\immediate\write\gtoutfile{--text follows this line--}
\immediate\write\gtoutfile{Proxy-for: \ifx\theasciiauthors\relax
\theauthors\else\theasciiauthors\fi\s<\ifx\theasciiemail\relax\theemail\else\theasciiemail\fi>}
\immediate\write\gtoutfile{\noexpand\\}
\immediate\write\gtoutfile{Authors: \ifx\theasciiauthors\relax
\theauthors\else\theasciiauthors\fi}
\immediate\write\gtoutfile{Title: \ifx\theasciititle\relax
\thetitle\else\theasciititle\fi}
\immediate\write\gtoutfile{Subj-class: GT or SG or MG etc}
\immediate\write\gtoutfile{MSC-class: \theprimaryclass\ifx\thesecondaryclass\relax\else, \thesecondaryclass\fi}
\immediate\write\gtoutfile{Journal-ref: Geom. Topol. \thevolumenumber
(\thevolumeyear) \startpage-\finishpage}
\immediate\write\gtoutfile{Comments: Published by Geometry and Topology at}
\immediate\write\gtoutfile{\s\s http://www.maths.warwick.ac.uk/gt/GTVol\thevolumenumber/paper\thepapernumber.abs.html}
\immediate\write\gtoutfile{\noexpand\\}
\immediate\write\gtoutfile{}
\ifx\theasciiabstract\relax
\immediate\write\gtoutfile{\theabstract}\else
\immediate\write\gtoutfile{\theasciiabstract}\fi
\immediate\write\gtoutfile{}
\immediate\write\gtoutfile{\noexpand\\}
\immediate\write\gtoutfile{}
\immediate\closeout\gtoutfile}}  

\def\maketitlepage{\maketitlep\makeheadfile}
\let\maketitle\maketitlepage

\def\ifplaintex{\expandafter\ifx\csname documentclass\endcsname\relax}

\ifplaintex 
\else
\fi

\expandafter\ifx\csname beginpicture\endcsname\relax
\expandafter\ifx\csname documentclass\endcsname\relax
\input pictex \else
\input prepictex \input pictex \input postpictex \fi\fi

\def\gt{{\mathsurround=0pt\it $\cal G\mskip-2mu$eometry \&\ 
$\cal T\!\!$opology}}        

\def\gtp{{\mathsurround=0pt\it $\cal G\mskip-2mu$eometry \&\ 
$\cal T\!\!$opology $\cal P\!$ublications}}  

\def\lognumber#1{\def\thelognumber{#1}}
\def\volumenumber#1{\def\thevolumenumber{#1}}
\def\papernumber#1{\def\thepapernumber{#1}}
\def\volumeyear#1{\def\thevolumeyear{#1}}

\def\pagenumbers#1#2{\def\startpage{#1}\def\finishpage{#2}}
\def\published#1{\def\publishdate{#1}}
\def\proposed#1{\def\theproposer{#1}}
\def\seconded#1{\def\theseconders{#1}}
\def\received#1{\def\receiveddate{#1}}

\def\accepted#1{\def\accepteddate{#1}}

\def\asciiaddress#1{\def\theasciiaddress{#1}}
\def\asciiemail#1{\def\theasciiemail{#1}}
\long\def\asciiabstract#1{\long\def\theasciiabstract{#1}}
\def\asciikeywords#1{\def\theasciikeywords{#1}}

\let\\\par\let\thelognumber\relax
\let\thevolumenumber\relax\let\thepapernumber\relax
\let\thevolumeyear\relax\let\thesamplenumber\relax\let\startpage\relax
\let\finishpage\relax\let\publishdate\relax\let\receiveddate\relax
\let\reviseddate\relax\let\accepteddate\relax\let\theasciititle\relax
\let\theasciiauthors\relax\let\theasciiaddress\relax
\let\theasciiabstract\relax\let\theasciikeywords\relax
\let\theasciiemail\relax\let\theshortauthors\relax\let\theshorttitle\relax

\long\def\maketitlep{   

\count0=\startpage

\gt\hfill      
\beginpicture
\setcoordinatesystem units <0.33truein, 0.33truein> point at 2.2 0.9
\setplotsymbol ({$\cal G$})
\plotsymbolspacing=9truept
\circulararc 315 degrees from 0 1 center at 0 0
\setplotsymbol ({$\cal T$})
\circulararc 315 degrees from 1 -1 center at 1 0
\endpicture
\break
{\small\ifx\thesamplenumber\relax 
Volume \else Sample
\fi\thevolumenumber\ (\thevolumeyear)
\startpage--\finishpage\nl
Published: \publishdate}
\vglue 0.5truein plus 0.4fil minus 0.1truein

{\parskip=0pt\leftskip 0pt plus 1fil\def\\{\par\smallskip}{\ifplaintex\large
\else\Large\fi\bf\thetitle}\par\medskip}   

\vglue 0pt plus 0.1fil 

{\parskip=0pt\leftskip 0pt plus 1fil\def\\{\par}{\sc\theauthors}
\par\medskip}

\vglue 0pt plus 0.1fil 

{\small\parskip=0pt\let\newline\\
{\leftskip 0pt plus 1fil\def\\{\par}{\sl\theaddress}\par}
\expandafter\ifx\theemail\relax    
\relax\else\vglue 5pt plus 0.02fil minus 2pt\def\\{\stdspace{\rm 
and}\stdspace} 
\cl{Email:\stdspace\tt\theemail}\fi
\ifx\theurl\relax                  
\relax\else\vglue 5pt plus 0.02fil minus 2pt\def\\{\stdspace{\rm 
and}\stdspace}
\cl{URL:\stdspace\tt\theurl}\fi\par}

\vglue 7pt plus 0.3fil minus 3pt

{\bf Abstract}
\vglue 5pt plus 0.1fil minus 2pt

\theabstract

\vglue 7pt plus 0.3fil minus 3pt

{\bf AMS Classification numbers}\quad Primary:\quad \theprimaryclass

Secondary:\quad \thesecondaryclass

\vglue 5pt plus 0.3fil minus 2pt

{\bf Keywords}\quad \thekeywords

\vglue 10pt plus 0.5fil minus 5pt

{\small  Proposed: \theproposer\hfill Received: \receiveddate\nl
Seconded: \theseconders\hfill 
\ifx\reviseddate\relax                         
Accepted: \accepteddate                        
\else
Revised: \reviseddate                          
\fi}
\eject
}       

\let\maketitlepage\maketitlep
\let\maketitle\maketitlepage

\font\phead=cmsl9 scaled 950
\font=cmsl9 scaled 1050
\font\pnum=cmbx10 scaled 913
\font\lnum=cmbx10 
\font\pfoot=cmsl9 scaled 950
\font=cmsl9 scaled 1050
\ifplaintex
\headline{\vbox to 0pt{\vskip -4.5mm\line{\small\phead\ifnum
\count0=\startpage ISSN 1364-0380 (on line)
1465-3060 (printed) \hfill {\pnum\folio}\else\ifodd\count0\def\\{ }%
\ifx\theshorttitle\relax\thetitle\else\theshorttitle\fi\hfill{\pnum\folio}
\else\def\\{ and }{\pnum\folio}\hfill\ifx\theshortauthors\relax\theauthors
\else\theshortauthors\fi\fi\fi}\vss}}
\footline{\vbox to 0pt{\vglue 0mm\line{\small\pfoot\ifnum\count0=\startpage
\copyright\ \gtp\hfill\else
\gt, Volume \thevolumenumber\ (\thevolumeyear)\hfill\fi}\vss
}}
\else
\makeatletter
\def\@oddhead{{\small\lhead\ifnum\count0=\startpage ISSN 1364-0380 (on line)
1465-3060 (printed) \hfill {\lnum\number\count0}\else\ifodd\count0
\def\\{ }\ifx\theshorttitle\relax \thetitle \else\theshorttitle\fi\hfill
{\lnum\number\count0}\else\def\\{ and }{\lnum\number\count0}
\hfill\ifx\theshortauthors\relax 
\theauthors\else\theshortauthors\fi\fi\fi}}\def\@evenhead{\@oddhead}
\def\@oddfoot{\small\lfoot\ifnum\count0=\startpage\copyright\ \gtp\hfill\else
\gt, Volume \thevolumenumber\ (\thevolumeyear)\hfill\fi}
\def\@evenfoot{\@oddfoot}
\makeatother
\fi

\newwrite\gtoutfile
\long\gdef\makeheadfile{  
{\def\\{, }\def\s{ }
\immediate\openout\gtoutfile head.xxx
\immediate\write\gtoutfile{To: math@arxiv.org}
\immediate\write\gtoutfile{Subject: put or rep NNNNN:pppp}
\immediate\write\gtoutfile{--text follows this line--}
\immediate\write\gtoutfile{Proxy-for: \ifx\theasciiauthors\relax
\theauthors\else\theasciiauthors\fi\s<\ifx\theasciiemail\relax\theemail\else\theasciiemail\fi>}
\immediate\write\gtoutfile{\noexpand\\}
\immediate\write\gtoutfile{Authors: \ifx\theasciiauthors\relax
\theauthors\else\theasciiauthors\fi}
\immediate\write\gtoutfile{Title: \ifx\theasciititle\relax
\thetitle\else\theasciititle\fi}
\immediate\write\gtoutfile{Subj-class: GT or SG or MG etc}
\immediate\write\gtoutfile{MSC-class: \theprimaryclass\ifx\thesecondaryclass\relax\else, \thesecondaryclass\fi}
\immediate\write\gtoutfile{Journal-ref: Geom. Topol. \thevolumenumber
(\thevolumeyear) \startpage-\finishpage}
\immediate\write\gtoutfile{Comments: Published by Geometry and Topology at}
\immediate\write\gtoutfile{\s\s http://www.maths.warwick.ac.uk/gt/GTVol\thevolumenumber/paper\thepapernumber.abs.html}
\immediate\write\gtoutfile{\noexpand\\}
\immediate\write\gtoutfile{}
\ifx\theasciiabstract\relax
\immediate\write\gtoutfile{\theabstract}\else
\immediate\write\gtoutfile{\theasciiabstract}\fi
\immediate\write\gtoutfile{}
\immediate\write\gtoutfile{\noexpand\\}
\immediate\write\gtoutfile{}
\immediate\closeout\gtoutfile}}  

\def\maketitlepage{\maketitlep\makeheadfile}
\let\maketitle\maketitlepage

\def\ifplaintex{\expandafter\ifx\csname documentclass\endcsname\relax}

\ifplaintex 
\else
\fi

\expandafter\ifx\csname beginpicture\endcsname\relax
\expandafter\ifx\csname documentclass\endcsname\relax
\input pictex \else
\input prepictex \input pictex \input postpictex \fi\fi

\def\gt{{\mathsurround=0pt\it $\cal G\mskip-2mu$eometry \&\ 
$\cal T\!\!$opology}}        

\def\gtp{{\mathsurround=0pt\it $\cal G\mskip-2mu$eometry \&\ 
$\cal T\!\!$opology $\cal P\!$ublications}}  

\def\lognumber#1{\def\thelognumber{#1}}
\def\volumenumber#1{\def\thevolumenumber{#1}}
\def\papernumber#1{\def\thepapernumber{#1}}
\def\volumeyear#1{\def\thevolumeyear{#1}}

\def\pagenumbers#1#2{\def\startpage{#1}\def\finishpage{#2}}
\def\published#1{\def\publishdate{#1}}
\def\proposed#1{\def\theproposer{#1}}
\def\seconded#1{\def\theseconders{#1}}
\def\received#1{\def\receiveddate{#1}}

\def\accepted#1{\def\accepteddate{#1}}

\def\asciiaddress#1{\def\theasciiaddress{#1}}
\def\asciiemail#1{\def\theasciiemail{#1}}
\long\def\asciiabstract#1{\long\def\theasciiabstract{#1}}
\def\asciikeywords#1{\def\theasciikeywords{#1}}

\let\\\par\let\thelognumber\relax
\let\thevolumenumber\relax\let\thepapernumber\relax
\let\thevolumeyear\relax\let\thesamplenumber\relax\let\startpage\relax
\let\finishpage\relax\let\publishdate\relax\let\receiveddate\relax
\let\reviseddate\relax\let\accepteddate\relax\let\theasciititle\relax
\let\theasciiauthors\relax\let\theasciiaddress\relax
\let\theasciiabstract\relax\let\theasciikeywords\relax
\let\theasciiemail\relax\let\theshortauthors\relax\let\theshorttitle\relax

\long\def\maketitlep{   

\count0=\startpage

\gt\hfill      
\beginpicture
\setcoordinatesystem units <0.33truein, 0.33truein> point at 2.2 0.9
\setplotsymbol ({$\cal G$})
\plotsymbolspacing=9truept
\circulararc 315 degrees from 0 1 center at 0 0
\setplotsymbol ({$\cal T$})
\circulararc 315 degrees from 1 -1 center at 1 0
\endpicture
\break
{\small\ifx\thesamplenumber\relax 
Volume \else Sample
\fi\thevolumenumber\ (\thevolumeyear)
\startpage--\finishpage\nl
Published: \publishdate}
\vglue 0.5truein plus 0.4fil minus 0.1truein

{\parskip=0pt\leftskip 0pt plus 1fil\def\\{\par\smallskip}{\ifplaintex\large
\else\Large\fi\bf\thetitle}\par\medskip}   

\vglue 0pt plus 0.1fil 

{\parskip=0pt\leftskip 0pt plus 1fil\def\\{\par}{\sc\theauthors}
\par\medskip}

\vglue 0pt plus 0.1fil 

{\small\parskip=0pt\let\newline\\
{\leftskip 0pt plus 1fil\def\\{\par}{\sl\theaddress}\par}
\expandafter\ifx\theemail\relax    
\relax\else\vglue 5pt plus 0.02fil minus 2pt\def\\{\stdspace{\rm 
and}\stdspace} 
\cl{Email:\stdspace\tt\theemail}\fi
\ifx\theurl\relax                  
\relax\else\vglue 5pt plus 0.02fil minus 2pt\def\\{\stdspace{\rm 
and}\stdspace}
\cl{URL:\stdspace\tt\theurl}\fi\par}

\vglue 7pt plus 0.3fil minus 3pt

{\bf Abstract}
\vglue 5pt plus 0.1fil minus 2pt

\theabstract

\vglue 7pt plus 0.3fil minus 3pt

{\bf AMS Classification numbers}\quad Primary:\quad \theprimaryclass

Secondary:\quad \thesecondaryclass

\vglue 5pt plus 0.3fil minus 2pt

{\bf Keywords}\quad \thekeywords

\vglue 10pt plus 0.5fil minus 5pt

{\small  Proposed: \theproposer\hfill Received: \receiveddate\nl
Seconded: \theseconders\hfill 
\ifx\reviseddate\relax                         
Accepted: \accepteddate                        
\else
Revised: \reviseddate                          
\fi}
\eject
}       

\let\maketitlepage\maketitlep
\let\maketitle\maketitlepage

\font\phead=cmsl9 scaled 950
\font=cmsl9 scaled 1050
\font\pnum=cmbx10 scaled 913
\font\lnum=cmbx10 
\font\pfoot=cmsl9 scaled 950
\font=cmsl9 scaled 1050
\ifplaintex
\headline{\vbox to 0pt{\vskip -4.5mm\line{\small\phead\ifnum
\count0=\startpage ISSN 1364-0380 (on line)
1465-3060 (printed) \hfill {\pnum\folio}\else\ifodd\count0\def\\{ }%
\ifx\theshorttitle\relax\thetitle\else\theshorttitle\fi\hfill{\pnum\folio}
\else\def\\{ and }{\pnum\folio}\hfill\ifx\theshortauthors\relax\theauthors
\else\theshortauthors\fi\fi\fi}\vss}}
\footline{\vbox to 0pt{\vglue 0mm\line{\small\pfoot\ifnum\count0=\startpage
\copyright\ \gtp\hfill\else
\gt, Volume \thevolumenumber\ (\thevolumeyear)\hfill\fi}\vss
}}
\else
\makeatletter
\def\@oddhead{{\small\lhead\ifnum\count0=\startpage ISSN 1364-0380 (on line)
1465-3060 (printed) \hfill {\lnum\number\count0}\else\ifodd\count0
\def\\{ }\ifx\theshorttitle\relax \thetitle \else\theshorttitle\fi\hfill
{\lnum\number\count0}\else\def\\{ and }{\lnum\number\count0}
\hfill\ifx\theshortauthors\relax 
\theauthors\else\theshortauthors\fi\fi\fi}}\def\@evenhead{\@oddhead}
\def\@oddfoot{\small\lfoot\ifnum\count0=\startpage\copyright\ \gtp\hfill\else
\gt, Volume \thevolumenumber\ (\thevolumeyear)\hfill\fi}
\def\@evenfoot{\@oddfoot}
\makeatother
\fi

\newwrite\gtoutfile
\long\gdef\makeheadfile{  
{\def\\{, }\def\s{ }
\immediate\openout\gtoutfile head.xxx
\immediate\write\gtoutfile{To: math@arxiv.org}
\immediate\write\gtoutfile{Subject: put or rep NNNNN:pppp}
\immediate\write\gtoutfile{--text follows this line--}
\immediate\write\gtoutfile{Proxy-for: \ifx\theasciiauthors\relax
\theauthors\else\theasciiauthors\fi\s<\ifx\theasciiemail\relax\theemail\else\theasciiemail\fi>}
\immediate\write\gtoutfile{\noexpand\\}
\immediate\write\gtoutfile{Authors: \ifx\theasciiauthors\relax
\theauthors\else\theasciiauthors\fi}
\immediate\write\gtoutfile{Title: \ifx\theasciititle\relax
\thetitle\else\theasciititle\fi}
\immediate\write\gtoutfile{Subj-class: GT or SG or MG etc}
\immediate\write\gtoutfile{MSC-class: \theprimaryclass\ifx\thesecondaryclass\relax\else, \thesecondaryclass\fi}
\immediate\write\gtoutfile{Journal-ref: Geom. Topol. \thevolumenumber
(\thevolumeyear) \startpage-\finishpage}
\immediate\write\gtoutfile{Comments: Published by Geometry and Topology at}
\immediate\write\gtoutfile{\s\s http://www.maths.warwick.ac.uk/gt/GTVol\thevolumenumber/paper\thepapernumber.abs.html}
\immediate\write\gtoutfile{\noexpand\\}
\immediate\write\gtoutfile{}
\ifx\theasciiabstract\relax
\immediate\write\gtoutfile{\theabstract}\else
\immediate\write\gtoutfile{\theasciiabstract}\fi
\immediate\write\gtoutfile{}
\immediate\write\gtoutfile{\noexpand\\}
\immediate\write\gtoutfile{}
\immediate\closeout\gtoutfile}}  

\def\maketitlepage{\maketitlep\makeheadfile}
\let\maketitle\maketitlepage

\def\ifplaintex{\expandafter\ifx\csname documentclass\endcsname\relax}

\ifplaintex 
\else
\fi

\expandafter\ifx\csname beginpicture\endcsname\relax
\expandafter\ifx\csname documentclass\endcsname\relax
\input pictex \else
\input prepictex \input pictex \input postpictex \fi\fi

\def\gt{{\mathsurround=0pt\it $\cal G\mskip-2mu$eometry \&\ 
$\cal T\!\!$opology}}        

\def\gtp{{\mathsurround=0pt\it $\cal G\mskip-2mu$eometry \&\ 
$\cal T\!\!$opology $\cal P\!$ublications}}  

\def\lognumber#1{\def\thelognumber{#1}}
\def\volumenumber#1{\def\thevolumenumber{#1}}
\def\papernumber#1{\def\thepapernumber{#1}}
\def\volumeyear#1{\def\thevolumeyear{#1}}

\def\pagenumbers#1#2{\def\startpage{#1}\def\finishpage{#2}}
\def\published#1{\def\publishdate{#1}}
\def\proposed#1{\def\theproposer{#1}}
\def\seconded#1{\def\theseconders{#1}}
\def\received#1{\def\receiveddate{#1}}

\def\accepted#1{\def\accepteddate{#1}}

\def\asciiaddress#1{\def\theasciiaddress{#1}}
\def\asciiemail#1{\def\theasciiemail{#1}}
\long\def\asciiabstract#1{\long\def\theasciiabstract{#1}}
\def\asciikeywords#1{\def\theasciikeywords{#1}}

\let\\\par\let\thelognumber\relax
\let\thevolumenumber\relax\let\thepapernumber\relax
\let\thevolumeyear\relax\let\thesamplenumber\relax\let\startpage\relax
\let\finishpage\relax\let\publishdate\relax\let\receiveddate\relax
\let\reviseddate\relax\let\accepteddate\relax\let\theasciititle\relax
\let\theasciiauthors\relax\let\theasciiaddress\relax
\let\theasciiabstract\relax\let\theasciikeywords\relax
\let\theasciiemail\relax\let\theshortauthors\relax\let\theshorttitle\relax

\long\def\maketitlep{   

\count0=\startpage

\gt\hfill      
\beginpicture
\setcoordinatesystem units <0.33truein, 0.33truein> point at 2.2 0.9
\setplotsymbol ({$\cal G$})
\plotsymbolspacing=9truept
\circulararc 315 degrees from 0 1 center at 0 0
\setplotsymbol ({$\cal T$})
\circulararc 315 degrees from 1 -1 center at 1 0
\endpicture
\break
{\small\ifx\thesamplenumber\relax 
Volume \else Sample
\fi\thevolumenumber\ (\thevolumeyear)
\startpage--\finishpage\nl
Published: \publishdate}
\vglue 0.5truein plus 0.4fil minus 0.1truein

{\parskip=0pt\leftskip 0pt plus 1fil\def\\{\par\smallskip}{\ifplaintex\large
\else\Large\fi\bf\thetitle}\par\medskip}   

\vglue 0pt plus 0.1fil 

{\parskip=0pt\leftskip 0pt plus 1fil\def\\{\par}{\sc\theauthors}
\par\medskip}

\vglue 0pt plus 0.1fil 

{\small\parskip=0pt\let\newline\\
{\leftskip 0pt plus 1fil\def\\{\par}{\sl\theaddress}\par}
\expandafter\ifx\theemail\relax    
\relax\else\vglue 5pt plus 0.02fil minus 2pt\def\\{\stdspace{\rm 
and}\stdspace} 
\cl{Email:\stdspace\tt\theemail}\fi
\ifx\theurl\relax                  
\relax\else\vglue 5pt plus 0.02fil minus 2pt\def\\{\stdspace{\rm 
and}\stdspace}
\cl{URL:\stdspace\tt\theurl}\fi\par}

\vglue 7pt plus 0.3fil minus 3pt

{\bf Abstract}
\vglue 5pt plus 0.1fil minus 2pt

\theabstract

\vglue 7pt plus 0.3fil minus 3pt

{\bf AMS Classification numbers}\quad Primary:\quad \theprimaryclass

Secondary:\quad \thesecondaryclass

\vglue 5pt plus 0.3fil minus 2pt

{\bf Keywords}\quad \thekeywords

\vglue 10pt plus 0.5fil minus 5pt

{\small  Proposed: \theproposer\hfill Received: \receiveddate\nl
Seconded: \theseconders\hfill 
\ifx\reviseddate\relax                         
Accepted: \accepteddate                        
\else
Revised: \reviseddate                          
\fi}
\eject
}       

\let\maketitlepage\maketitlep
\let\maketitle\maketitlepage

\font\phead=cmsl9 scaled 950
\font=cmsl9 scaled 1050
\font\pnum=cmbx10 scaled 913
\font\lnum=cmbx10 
\font\pfoot=cmsl9 scaled 950
\font=cmsl9 scaled 1050
\ifplaintex
\headline{\vbox to 0pt{\vskip -4.5mm\line{\small\phead\ifnum
\count0=\startpage ISSN 1364-0380 (on line)
1465-3060 (printed) \hfill {\pnum\folio}\else\ifodd\count0\def\\{ }%
\ifx\theshorttitle\relax\thetitle\else\theshorttitle\fi\hfill{\pnum\folio}
\else\def\\{ and }{\pnum\folio}\hfill\ifx\theshortauthors\relax\theauthors
\else\theshortauthors\fi\fi\fi}\vss}}
\footline{\vbox to 0pt{\vglue 0mm\line{\small\pfoot\ifnum\count0=\startpage
\copyright\ \gtp\hfill\else
\gt, Volume \thevolumenumber\ (\thevolumeyear)\hfill\fi}\vss
}}
\else
\makeatletter
\def\@oddhead{{\small\lhead\ifnum\count0=\startpage ISSN 1364-0380 (on line)
1465-3060 (printed) \hfill {\lnum\number\count0}\else\ifodd\count0
\def\\{ }\ifx\theshorttitle\relax \thetitle \else\theshorttitle\fi\hfill
{\lnum\number\count0}\else\def\\{ and }{\lnum\number\count0}
\hfill\ifx\theshortauthors\relax 
\theauthors\else\theshortauthors\fi\fi\fi}}\def\@evenhead{\@oddhead}
\def\@oddfoot{\small\lfoot\ifnum\count0=\startpage\copyright\ \gtp\hfill\else
\gt, Volume \thevolumenumber\ (\thevolumeyear)\hfill\fi}
\def\@evenfoot{\@oddfoot}
\makeatother
\fi

\newwrite\gtoutfile
\long\gdef\makeheadfile{  
{\def\\{, }\def\s{ }
\immediate\openout\gtoutfile head.xxx
\immediate\write\gtoutfile{To: math@arxiv.org}
\immediate\write\gtoutfile{Subject: put or rep NNNNN:pppp}
\immediate\write\gtoutfile{--text follows this line--}
\immediate\write\gtoutfile{Proxy-for: \ifx\theasciiauthors\relax
\theauthors\else\theasciiauthors\fi\s<\ifx\theasciiemail\relax\theemail\else\theasciiemail\fi>}
\immediate\write\gtoutfile{\noexpand\\}
\immediate\write\gtoutfile{Authors: \ifx\theasciiauthors\relax
\theauthors\else\theasciiauthors\fi}
\immediate\write\gtoutfile{Title: \ifx\theasciititle\relax
\thetitle\else\theasciititle\fi}
\immediate\write\gtoutfile{Subj-class: GT or SG or MG etc}
\immediate\write\gtoutfile{MSC-class: \theprimaryclass\ifx\thesecondaryclass\relax\else, \thesecondaryclass\fi}
\immediate\write\gtoutfile{Journal-ref: Geom. Topol. \thevolumenumber
(\thevolumeyear) \startpage-\finishpage}
\immediate\write\gtoutfile{Comments: Published by Geometry and Topology at}
\immediate\write\gtoutfile{\s\s http://www.maths.warwick.ac.uk/gt/GTVol\thevolumenumber/paper\thepapernumber.abs.html}
\immediate\write\gtoutfile{\noexpand\\}
\immediate\write\gtoutfile{}
\ifx\theasciiabstract\relax
\immediate\write\gtoutfile{\theabstract}\else
\immediate\write\gtoutfile{\theasciiabstract}\fi
\immediate\write\gtoutfile{}
\immediate\write\gtoutfile{\noexpand\\}
\immediate\write\gtoutfile{}
\immediate\closeout\gtoutfile}}  

\def\maketitlepage{\maketitlep\makeheadfile}
\let\maketitle\maketitlepage

\lognumber{143}
\volumenumber{5}\papernumber{9}\volumeyear{2001}
\pagenumbers{287}{318}
\accepted{20 March 2001}
\received{13 October 2000}

\published{24 March 2001\nl{\sl Version 2 published 8 June 2002:\nl
 Corrections to equation (2) page 295, to the first equation\nl in
Proposition 2.1 and to the tables on page 318}}

\proposed{Robion Kirby}
\seconded{Yasha Eliashberg, Simon Donaldson}

\usepackage{amsmath,amsthm,amsfonts}

\newtheorem{thm}{Theorem}[section]
\newtheorem{theorem}[thm]{Theorem}

\newtheorem{lem}[thm]{Lemma}
\newtheorem{lemma}[thm]{Lemma}

\newtheorem{proposition}[thm]{Proposition}

\newtheorem{conjecture}[thm]{Conjecture}

\newcounter{Pcounter}

\newcounter{Ccounter}

\newtheorem{Plemmas}[Pcounter]{Lemma}
\newtheorem{Clemmas}[Ccounter]{Lemma}

\theoremstyle{remark}
\newtheorem{remark}[thm]{Remark}
\newtheorem{defn}[thm]{Definition}
\newtheorem{definition}[thm]{Definition}

\newcommand{\cnums} {{\mathbf C}}          
\newcommand{\nnums} {{\mathbf N}}		
\newcommand{\rnums} {{\mathbf R}}		
\newcommand{\znums} {{\mathbf Z}}		
\newcommand{\qnums} {{\mathbf Q}}		

\newcommand{\Hom}{\operatorname{Hom}}

\newcommand{\im}{\operatorname{Im}}
\newcommand{\Coeff}[1]{\operatorname{Coeff}_{#1}}

\newcommand{\ie}{ie,}

\newcommand{\til}[1]{\widetilde{#1}}

\renewcommand{\P}{\mathbf{P}}
\newcommand{\M}{\overline{{M}}}

\newcommand{\et}{\acute{e}t}
\renewcommand{\O}{\mathcal{O}}
\renewcommand{\arcsin}{\operatorname{arcsin}}
\renewcommand{\hat}{\widehat}
\newcommand{\Sym}{\operatorname{Sym}}

\begin{document}

\title{BPS states of curves in Calabi--Yau 3--folds}
\author{Jim Bryan\\Rahul Pandharipande}
\address{
Department of Mathematics,Tulane University\\
6823 St. Charles Ave, New Orleans, LA 70118, USA}
\secondaddress{Department of Mathematics, California Institute of Technology\\
Pasadena, CA 91125, USA}
\asciiaddress{
Department of Mathematics,Tulane University\\
6823 St. Charles Ave, New Orleans, LA 70118, USA\\
Department of Mathematics, California Institute of Technology\\
Pasadena, CA 91125, USA}
\email{jbryan@math.tulane.edu}
\secondemail{rahulp@its.caltech.edu}
\asciiemail{jbryan@math.tulane.edu, rahulp@its.caltech.edu}

\begin{abstract}
The Gopakumar--Vafa integrality conjecture is defined and studied for the
local geometry of a super-rigid curve in a Calabi--Yau 3--fold.  The
integrality predicted in Gromov--Witten theory by the Gopakumar--Vafa BPS
count is verified in a natural series of cases in this local geometry.  The
method involves Gromov--Witten computations, M\"obius inversion, and a
combinatorial analysis of the numbers of \'etale covers of a curve.
\end{abstract}

\asciiabstract{The Gopakumar-Vafa conjecture is defined and studied
for the local geometry of a curve in a Calabi-Yau 3-fold.  The
integrality predicted in Gromov-Witten theory by the Gopakumar-Vafa
BPS count is verified in a natural series of cases in this local
geometry. The method involves Gromov-Witten computations, Mobius
inversion, and a combinatorial analysis of the numbers of etale covers
of a curve.}

\asciikeywords{Gromov-Witten invariants, BPS states, Calabi-Yau 3-folds}

\maketitlepage

\section{Introduction and Results}

\subsection{Gromov--Witten and BPS invariants}\label{subsec: GW and BPS
invs} Let $X$ be a Calabi--Yau 3--fold and let $N_{\beta }^{g} (X)$ be the
0--point genus $g$ Gromov--Witten invariant of $X$ in the curve
class $\beta \in
H_{2} (X,\znums )$. From considerations in M--theory, Gopakumar and
Vafa express the invariants $N_{\beta }^{g} (X)$ in terms of integer
invariants $n_{\beta }^{g} (X)$ obtained by BPS state counts \cite{Go-Va}.
The Gopakumar--Vafa formula may be viewed as providing a
definition of the BPS state counts $n_{\beta}^g (X)$ in terms of the
Gromov--Witten invariants.

\begin{definition}\label{defn: BPS invariantss via G-V formula}
Define the \emph{Gopakumar--Vafa BPS invariants} 
$n_{\beta }^{{r}} (X)$ by the formula: 
\begin{equation}
\label{zzzz}
\sum _{\beta \neq 0}\sum _{{g}\geq 0} N_{\beta
}^{{g}} (X)t^{2{g}-2}q^{\beta }=\sum _{\beta \neq 0}\sum _{{g}\geq 0}
n^{{g}}_{\beta } (X)\sum _{k>0}\tfrac{1}{k} \left(2 \sin (
\tfrac{kt}{2}) \right)^{2{g}-2}q^{k\beta }.
\end{equation}
Matching the coefficients of the two series yields equations determining
$n_{\beta }^{g} (X)$ recursively in terms of $N_{\beta }^{g} (X)$ (see
Proposition~\ref{prop: general inversion of GV} for an explicit inversion
of this formula).
\end{definition}

From the above definition, there is no (mathematical) reason to expect
$n_{\beta }^{{g}} (X)$ to be an integer. Thus, the physics makes the
following prediction.

\begin{conjecture}\label{conj: BPS invariants are integers}
The BPS invariants are integers:
\[
n_{\beta }^{{g}} (X)\in \znums .
\]
Moreover, for any fixed $\beta $, $n_{\beta }^{g} (X)=0$ for $g>\!\!>0$.
\end{conjecture}

\begin{remark}\label{rem: progress on defn of BPS invariants}
By the physical arguments of Gopakumar and Vafa, the BPS invariants should
be directly defined via the cohomology of the $D$--brane moduli space.
First, the $D$--brane moduli space $\hat{M}$ should be defined with a
natural morphism $\hat{M}\to M$ to a moduli space $M$ of curves in $X$ in
the class $\beta$.  The fiber of $\hat{M}\to M$ over each curve $C\in M$
should parameterize flat line bundles on $C$. Furthermore, there should
exist an $\mathfrak{sl_{2}\oplus sl_{2}}$ representation on $H^{*}
(\hat{M},\cnums )$ such that the diagonal and right actions are the usual
$\mathfrak{sl_{2}}$ Lefschetz representations on $H^*(\hat{M},\cnums)$ and
$H^*(M,\cnums)$ respectively --- assuming $\hat{M}$ and $M$ are compact,
nonsingular, and K\"ahler.  The BPS state counts $n_{\beta }^{g} (X)$ are
then the coefficients in the decomposition of the left (fiberwise)
$\mathfrak{sl_{2}}$ representation $H^{*} (\hat{M},\cnums)$ in the basis
given by the cohomologies of the algebraic tori.  After these foundations
are developed, Equation~(\ref{zzzz}) should be {\em proven} as the basic
result relating Gromov--Witten theory to the BPS invariants.

The correct mathematical definition of the D--brane moduli space is unknown
at present, although there has been recent progress in case the curves move
in a surface $S\subset X$ (see \cite{HST}, \cite{KKV},
\cite{Klemm-Zaslow}). The nature of the D--brane moduli space in the case
where there are non-reduced curves in the family $M$ is not well
understood. The fiber of $\hat{M}\to M$ over a point corresponding to a
non-reduced curve may involve higher rank bundles on the reduction of the
curve.  It has been recently suggested by Hosono, Saito, and Takahashi
\cite{HST2} that the $\mathfrak{sl_{2}\oplus sl_{2}}$ representation can be
constructed in general via intersection cohomology and the
Beilinson--Bernstein--Deligne spectral sequence \cite{Be-Be-De}.
\end{remark}

\begin{remark}
An extension of formula (\ref{zzzz}) conjecturally defining integer
invariants for arbitrary 3--folds (not necessarily Calabi--Yau) has been
found in \cite{Pa}, \cite{Pandharipande-Vakil}.  Some predictions in the
non Calabi--Yau case have been verified in \cite{Br-E1test}.  Though it is
not yet known how the relevant physical arguments apply to the non
Calabi--Yau geometries, one may hope a mathematical development will provide
a unified approach to all 3--folds.
\end{remark}

The physical discussion suggests that the BPS invariants will be a sum of
integer contributions coming from each component of the D--brane moduli
space (whatever space that may be).
 One obvious source of such components occurs when
the curves parameterized by $M$ are rigid or lie in a fixed surface.  The
moduli space of stable maps has corresponding components given by those
maps whose image is the rigid curve or respectively lies in the fixed
surface. These give rise to the notion of ``local Gromov--Witten
invariants'' and we expect that the corresponding ``local BPS invariants''
will be integers.

\subsection{Local contributions}\label{subsec: local contributions} In this
paper we are interested in the contributions of an isolated curve $C\subset
X$ to the Gromov--Witten invariants $N_{d[C]}^{g} (X)$ and the BPS
invariants $n_{d[C]}^{g} (X)$.

To discuss the local contributions of a curve (also often called ``multiple
cover contributions''), we make the following definitions:
\begin{definition}\label{defn: local inv's, super-rigidity}
Let $C\subset X$ be a curve and let $M_{C}\subset \M _{g} (X,d[C])$ be the
locus of maps whose image is $C$. Suppose that   $M_{C}$ is an open
component of $\M _{g} (X,d[C])$.  Define the \emph{local 
Gromov--Witten invariant}, 
$N_{d}^{g} (C\subset X)\in
\qnums $ by the evalution of the well-defined restriction of
 $[\M _{g}
(X,d[C])]^{vir}$ to $H_{0} (M_{C},\qnums )$.  
\end{definition}

\begin{definition}
Let $C\subset X$ satisfy the conditions of Definition \ref
{defn: local inv's, super-rigidity}. 
If  $$M_{C}\cong \M _{g} (C,d)$$ then $C$ is said to be
\emph{$(d,g)$--rigid}. If $C$ is $(d,g)$--rigid for all $d$ and $g$, then
$C$ is \emph{super-rigid}.
\end{definition}

For example, a nonsingular rational curve with normal bundle $\O (-1)\oplus \O
(-1)$ is super-rigid. An elliptic curve $E\subset X$ is super-rigid if and
only if $N_{E/X}\cong L\oplus L^{-1}$ where $L\to E$ is a flat line bundle
such that no power of $L$ is trivial (see \cite{Pa}). 
An example where $M_{C}$ is an open
component but $M_{C}\not\cong \M _{g} (C,d)$ is the case where $C\subset X$
is a contractable, smoothly embedded $\cnums \P ^{1}$ with $N_{C/X}\cong \O
\oplus \O (-2)$. In this case $M_{C}$ has non-reduced structure coming from
the (obstructed) infinitesimal deformations of $C$ in the $\O $ direction
of $N_{C/X}$ (see \cite{BKL} for the computation of $N^{g}_{d} (C\subset
X)$ in this case).

The existence of genus $g$ curves in $X$ with $(d,g+h)$--rigidity is likely
to be a subtle question in the algebraic geometry of Calabi--Yau 3--folds.
On the other hand, these rigidity issues may be less delicate in the
symplectic setting. For a generic almost complex structure on $X$, it is
reasonable to hope super-rigidity will hold for any pseudo-holomorphic
curve in $X$.

Let $h\geq 0$ and suppose a nonsingular genus $g$ curve $C_{g}\subset X$ is
$(d,g+h)$--rigid. Then $N_{d}^{g+h} (C_{g}\subset X)$ can be expressed as
the integral of an Euler class of a bundle over $[\M _{g+h}
(C_{g},d)]^{vir}$. Let $\pi \co U\to \M _{g+h} (C_{g},d)$ be the universal
curve and let $f\co U\to C_{g}$ be the universal map. Then
\[
N_{d}^{g+h} (C_{g}\subset X)\cong \int _{[\M _{g+h}
(C_{g},d)]^{vir}} c (R^{1}\pi _{*}f^{*} (N_{C/X})).
\]
In fact, we can rewrite the above integral in the
following form: 
\begin{align*}
\int c (R^{1}\pi _{*}f^{*}N_{C/X}) 
&=\int c (R^{\bullet }\pi _{*}f^{*}N_{C/X}[1])\\
&=\int c (R^{\bullet }\pi _{*}f^{*} (\O _{C}\oplus \omega _{C})[1])
\end{align*}
where all the integrals are over $[\M _{g+h} (C_{g},d)]^{vir}$. The first
equality holds because $(d,g+h)$--rigidity implies that $R^{0}\pi
_{*}f^{*}N_{C/X}$ is 0. The second equality holds because $N_{C/X}$ deforms
to $\O _{C}\oplus \omega _{C}$, the sum of the trivial sheaf and the
canonical sheaf (this follows from an easily generalization of the argument
at the top of page~497 in \cite{Pa}). The last integral depends only upon
$g$, $h$, and $d$. We regard this formula as defining the idealized
multiple cover contribution of a genus $g$ curve by maps of degree $d$ and
genus $g+h$.

We will denote this idealized contribution by the following notation:
\[
N_{d}^{h} (g):=\int _{[\M _{g+h} (C_{g},d)]^{vir}}c (R^{\bullet }\pi
_{*}f^{*} (\O _{C}\oplus \omega _{C})[1]).
\]
From the previous discussion, $N_{d}^{h} (g)=N_{d}^{g+h} (C_{g})$ for any
nonsingular, $(d,g+h)$--rigid, genus $g$ curve $C_{g}$.

We define the local BPS invariants in terms of the local Gromov--Witten
invariants via the Gopakumar--Vafa formula.

\begin{defn}\label{defn: local BPS invariants}
Define the \emph{local BPS invariants} $n_{d}^{h} (g)$ in terms of
the local Gromov--Witten invariants by the formula
\[
\sum _{\beta \neq 0}\sum _{h\geq 0} N_{d
}^{h} (g)t^{2 (g+h-1)}q^{d }=\sum _{d\neq 0}\sum _{h\geq 0}
n^{h}_{d } (g)\sum _{k>0}\tfrac{1}{k} \left(2 \sin (
\tfrac{kt}{2}) \right)^{2 (g+h-1)}q^{kd }.
\]
\end{defn}

The local Gromov--Witten invariants $N_{d}^{h} (g)$ are in general difficult
to compute. For $g=0$, these integrals were computed in \cite{Fa-Pa}. In
terms of local BPS invariants, these calculations yield:
\begin{align*}
n_{d}^{h} (0)&=
\begin{cases}1&\text{for $d=1$ and $h=0$,}\\0&\text{otherwise.} \end{cases}\\
\end{align*}
For $g=1$, complete results have also been obtained \cite{Pa}:
\begin{align*}
n_{d}^{h} (1)&=
\begin{cases}1&\text{for $d\geq 1$ and $h=0$,}\\0&\text{otherwise.} 
\end{cases}\\
\end{align*}
The local invariants of a super-rigid nodal rational curve as well as the
local invariants of contractable (non-generic) embedded rational curves were
determined in \cite{BKL}.

In this paper we compute certain contributions to the local Gromov--Witten
invariants $N_{d}^{h} (g)$ for $g>1$ and we determine the corresponding
contributions to the BPS invariants $n^{h}_{d} (g)$. We prove the
integrality of these contributions. In the appendix, we provide tables
giving explicit values for $n_{d}^{h} (g)$.

\subsection{Results}\label{subsec: results}
The contributions to $N_{d}^{h} (g)$ we compute are those that come from
maps $[f:$\break$D \rightarrow C]$ satisfying either of following conditions:
\begin{enumerate}
\item[(i)] A single component of the domain is an \'etale cover of $C$
(with any number of auxiliary collapsed components {\em simply} attached to the
\'etale com\-ponent). 
\item[(ii)] The map $f$ has exactly two branch points (and no collapsed
components).
\end{enumerate}

The type (i) contributions, the {\em \'etale invariants}, correspond to the
first level in a natural grading on the set of local Gromov--Witten
invariants which will be discussed in Section~\ref{subsec: primitive GW
inv's}.  We use an elementary observation to reduce the computation of the
\'etale invariants to the computation of the degree 1 local invariants, \ie
$n^{h}_{1} (g)$. The computation of the degree 1 invariants was done
previously by the second author in \cite{Pa}. The observation that we use,
while elementary, seems useful enough to formalize in a general
setting. This we do by the introduction of \emph{primitive Gromov--Witten
invariants} in Section~\ref{subsec: primitive GW inv's}.


The type (ii) contribution we compute by a
Grothendieck--Riemann--Roch calculation which is carried out in
Section~\ref{sec: the GRR computation}. 

\subsubsection{Type (i) contributions (\'etale contributions)}

\begin{definition}\label{defn: etale invariants}
We define $\M _{g+h}^{\et } (C,d)\subset \M _{g+h} (C,d)$ to be the
union of the moduli components  corresponding to 
stable maps $\pi\co  D \rightarrow C$  satisfying:
\begin{enumerate} \item[(a)] $D$ contains a unique component $C'$ of degree
$d$, \'etale over $C$, while all other components are degree 0.
\item[(b)]  All
$\pi$--collapsed components are all simply attached to $C'$ (the vertex in
the dual graph of the domain curve corresponding to $C'$ does not contain a
cycle). 
\end{enumerate}
We define the \'etale Gromov--Witten invariants by
\[
N_{d}^{h} (g)^{\et}:=\int _{[\M _{g+h} (C,d)^{\et}]^{vir}} 
c
(R^{\bullet }\pi _{*}f^{*} (\mathcal{O}_{C}\oplus \omega _{C})[1])
\]
and we define the \'etale BPS invariants $n_{d}^{h} (g)^{\et}$ in terms of
$N_{d}^{h} (g)^{\et}$ via the Gopakumar--Vafa formula as before.
\end{definition}

As we will explain in Section~\ref{sec: inversion of GV formula}, any
Gromov--Witten invariant can be written in terms of \emph{primitive}
Gromov--Witten invariants. The \'etale invariants exactly correspond to
those that can be expressed in terms of degree 1 primitive invariants.

Our main two Theorems concerning the \'etale BPS invariants give an explicit
formula for $n_{d}^{h} (g)^{\et}$ and prove they are integers.

\begin{thm}\label{thm: formula for etale BPS invariants}
Let $C_{n,g}$ be number of degree $n$, connected, complete, \'etale covers
of a curve of genus $g$, each counted by the reciprocal of the number of
automorphisms of the cover. Let $\mu $ be the M\"obius function: $\mu (n)=
(-1)^{a}$ where $a$ is the number of prime factors of $n$ if $n$ is
square-free and $\mu (n)=0$ if $n$ is not square-free.  Then the \'etale
BPS invariants are given as the coefficients of the following polynomial:
\[
\sum _{h\geq 0}n_{d}^{h} (g)^{\et} y^{h+g-1}=
\sum _{k|d} k\mu
(k)\; C_{\frac{d}{k},g }\;   P_{k} (y) ^{\frac{d (g-1)}{k}}
\]
where the polynomial $P_{k} (y)$ is defined\footnote{\textbf{Warning:} This
definition of $P_{k} (y)$ differs from the one in \cite{Br-RIMS} by a
factor of $y$.} by
\[
P_{k} (4\sin^{2}t)=4\sin ^{2} (kt)
\] 
which by Lemma~\ref{prop: formula for the poly P} is given explicitly by
\[
P_{k} (y)=\sum _{a=1}^{k}-\frac{k}{a}\binom{a+k-1}{2a-1} (-y)^{a}.
\]
\end{thm}

\begin{theorem}\label{thm: etale BPS invariants are integers}
The \'etale BPS invariants are integers: $n_{d}^{h} (g)^{\et}\in
\znums $. 
\end{theorem}

We note that $C_{n,g}$ is not integral in general, for example $C_{2,g}=
(2^{2g}-1)/2$.  We also note that the formula given by the Theorem
\ref{thm: formula for etale BPS invariants} shows that for fixed $d$ and
$g$, $n_{d}^{h} (g)^{\et}$ is non-zero only if $0\leq h\leq (d-1)
(g-1)$. See Table 1 for explicit values of $n_{d}^{h} (g)^{\et }$ for small
$d$, $g$, and $h$.

There is a range where the \'etale contributions are the only contributions
to the full local BPS invariant. 
\begin{lemma}\label{lem: range where etale BPS inv equal usual}
Let $d_{min}$ be the smallest divisor $d'$ of $d$ that is not 1 and such
that $\mu (\frac{d}{d'})\neq 0$, then
\begin{equation*}
n_{d}^{h} (g)=n^{h}_{d} (g)^{\et }\quad \text{for all $h\leq (d_{min}-1)
(g-1)$.}
\end{equation*}
\end{lemma}
\proof This follows from Equation~\ref{eqn: inverted GV, formula
for each coef} (in Section~\ref{sec: inversion of GV formula}) and the
simple geometric fact that a degree $d$ stable map $f\co D_{g+h}\to C_{g}$
must be of type (i) if $h\leq (d-1) (g-1)$ or if $d=1$.\qed

\begin{remark}\label{rem: D-bran moduli space should have ``etale'' component}
\emph{A priori} there is no reason (even physically) to expect that the
\'etale invariants $n^{h}_{d} (g)^{\et }$ are integers outside of the range
where $n_{d}^{h} (g)^{\et }=n_{d}^{h} (g)$. Theorem~\ref{thm: etale BPS
invariants are integers} is very suggestive that the D--brane moduli space
has a distinguished component (or components) corresponding to these
\'etale contributions. Furthermore, our results suggest that this component
has dimension $d (g-1)+1$ and has a product decomposition (at least
cohomologically) with one factor a complex torus of dimension $g$.
\end{remark}

Theorem \ref{thm: formula for etale BPS invariants} follows from the
computation of $N^{h}_{d} (g)^{\et}$ by a (reasonably straightforward)
inversion of the Gopakumar--Vafa formula that is carried out in
Section~\ref{sec: inversion of GV formula}.  Theorem \ref{thm: etale BPS
invariants are integers} is proved directly from the formula given in
Theorem \ref{thm: formula for etale BPS invariants} and turns out to be
rather involved. It depends on somewhat delicate congruence properties of
the polynomials $P_{l} (y)$ and the number of covers $C_{d,g}$. These are
proved in Section~\ref{sec: pf of integrality}.

\subsubsection{Type (ii) contributions}\label{subsec: type ii contributions}

There is another situation where $\M _{g+h} (C_{g},d)$ has a distinguished
open component. If
\[
h= (d-1) (g-1)+1,
\]
then there are exactly two open components, namely the \'etale component
$M^{\et }$ and one other $\til{M}\subset \M _{g+h} (C_{g},d)$.  The
generic points of $\til{M}$ correspond to maps of nonsingular curves with
exactly two simple ramification points. Let $\til{N}_{g} (d)$ be the
corresponding contribution to the Gromov--Witten invariants so that
\[
N_{d}^{(d-1) (g-1)+1} (g) = N_{d}^{(d-1) (g-1)+1} (g)^{\et } + \til{N}_{d} (g).
\]
The component $\til{M}$ admits a finite morphism to $\operatorname{Sym}^{2}
(C_{g})$ given by sending a map to its branched locus (see
\cite{Fa-Pa} for the existence of such a morphism). 

We compute the invariant $\til{N}_{d} (g)$ in Section~\ref{sec: the GRR
computation} by a Grothendieck--Riemann--Roch (GRR) computation. The relative
Todd class required by GRR is computed using the formula of Mumford
\cite{Mumford} adapted to the context of stable maps (see \cite{Fa-Pa}
Section~1.1). The intersections in the GRR formula are computed by pushing
forward to $\operatorname{Sym^{2}} (C_{g})$. The result of this computation
(which is carried out in Section~\ref{sec: the GRR computation}) is the
following:

\begin{thm}\label{thm: GRR computation}
$$\til{N}_{d} (g)=\int _{\til{M}}c (R^{\bullet }\pi _{*}f^{*}
(\mathcal{O}_{C_{g}}\oplus \omega _{C_{g}})[1]) =
\frac{g-1}{8} \Big( (g-1)D_{d,g} - D_{d,g}^* -\frac{1}{27}
D_{d,g}^{**}\Big).$$
\end{thm}

The numbers
$D_{d,g}$, $D_{d,g}^*$, and
$D_{d,g}^{**}$ are the following Hurwitz numbers of covers of the
curve $C_g$.

\begin{enumerate}
\item[$\bullet$]$D_{d,g}$  is the number of connected, degree $d$
covers of $C_{g}$ simply  branched over 2 distinct fixed points of $C_{g}$.
\item[$\bullet$]
$D^*_{d,g}$ is the number of connected, degree $d$, covers of $C_{g}$
with  1 node lying over a fixed point of $C_{g}$. 
\item[$\bullet$]
$D^{**}_{d,g}$ is the number of connected, degree $d$
covers of $C_{g}$ with  1 double ramification point over a fixed point
of $C_{g}$.
\end{enumerate}

The covers are understood to be \'etale away from the imposed
ramification. Also, $D_{d,g}$, $D^*_{d,g}$, and $D^{**}_{d,g}$ are all
counts weighted by the reciprocal of the number of automorphisms of the
covers.

There is an additional Hurwitz number $D^{***}_{g,d}$
which is natural to consider together with the three above:
\begin{enumerate}
\item[$\bullet$] $D^{***}_{d,g}$ is the number 
of connected, degree $d$ covers of $C_{g}$ with
2 distinct ramification points in the domain lying
over a fixed  point of $C_{g}$. 
\end{enumerate}
However, $D_{d,g}^{***}$
is determined from the previous Hurwitz numbers by the degeneration
relation:
\begin{equation}\label{eqn: Hurwitz degeneration relation}
D_{d,g}= D_{d,g}^* + 3D_{d,g}^{**} + 2D_{d,g}^{***}            
\end{equation}
(see \cite{Harris-Morrison}).  Theorem~\ref{thm: GRR computation} therefore
involves all of the independent covering numbers which appear in this 2
branch point geometry (see Table~3 for some explicit values of these
numbers).

Theorem~\ref{thm: GRR computation} can be used to extend the range where we
can compute the full local BPS invariants.  Lemma~\ref{lem: range where
etale BPS inv equal usual} generalizes to
\begin{lemma}\label{lem: local BPS invs in known range (including GRR)}
Let $d_{min} $ be defined as in Lemma~\ref{lem: range where etale BPS inv
equal usual}, then
\begin{equation*}
n_{d}^{h} (g)=\begin{cases}
n_{d}^{h} (g)^{\et }&\text{for all $h\leq (d_{min}-1) (g-1)$}\\
n_{d}^{h} (g)^{\et }+\epsilon \til{N}_{d_{min}} (g)&\text{for $h=
(d_{min}-1) (g-1)+1$}
\end{cases}
\end{equation*}
where $\epsilon $ is the rational number given by Equation~\ref{eqn: inverted GV, formula for each coef}, \ie 
\[
\epsilon =\mu (\tfrac{d}{d_{min}}) (\tfrac{d}{d_{min}})^{d_{min} (g-1)+2}.
\]
\end{lemma}

For example, if $d$ is prime, then $d_{min}=d$ and $\epsilon =1$. See
Table~2 for explicit values of $n_{d}^{h} (g)$ for small $d$, $g$, and
$h$.

Since $n^{h}_{d} (g)^{\et }\in \znums $ by Theorem~\ref{thm: etale BPS
invariants are integers}, the integrality conjecture predicts that
$\epsilon \til{N}_{d} (g)\in \znums $. In light of our formula in
Theorem~\ref{thm: GRR computation}, this leads to congruences that are
conjecturally satisfied by the Hurwitz numbers $D_{d,g}$, $D_{d,g}^{*}$,
and $D_{d,g}^{**}$.

\begin{conjecture}\label{conj: congruence of hurwitz nums}
Let $\Upsilon _{d,g}=216\til{N}_{d} (g)$, that is
\[
\Upsilon _{d,g}= (g-1)\left(27 (g-1)D_{d,g}-27D^{*}_{d,g}-D^{**}_{d,g}
\right).
\]
Suppose that $d$ is not divisible by 4, 6, or 9. Then, 
\[
\Upsilon_{d,g}   \equiv 0 \pmod {216}.
\]
\end{conjecture}

Although $D_{d,g}$, $D_{d,g}^{*}$, and $D_{d,g}^{**}$ are not \emph{a
priori} integers, it is proven in \cite{Br-RIMS} that $\Upsilon _{d,g}\in
\znums $. It is also proven in \cite{Br-RIMS} that Conjecture~\ref{conj:
congruence of hurwitz nums} holds for $d=2$ and $d=3$.

\begin{remark}
Various congruence properties of $C_{d,g}$ (the number of degree $d$
connected \'etale covers) will also be used in the proof of the integrality of
the \'etale BPS invariants $n_{d}^{h} (g)^{\et }$ (see Lemma~\ref{lem: mod
p properties of the c's}). We speculate that these and the above conjecture
are the beginning of a series of congruence properties of general Hurwitz
numbers that are encoded in the integrality of the local BPS invariants.
\end{remark}

\subsection{Acknowledgements} The research presented here began during a
visit to the ICTP in Trieste in summer of 1999. We thank M Aschbacher,
C Faber, S Katz, V Moll, C Vafa, R Vakil, and E Zaslow for many
helpful discussions.  The authors were supported by Alfred P Sloan
Research Fellowships and NSF grants DMS-9802612, DMS-9801574, and
DMS-0072492.

\section{Inversion of the Gopakumar--Vafa formula}\label{sec: inversion of
GV formula}

In this section we invert the Gopakumar--Vafa formula in general to give an
explicit expression for the BPS invariants in terms of the Gromov--Witten
invariants. We then introduce the notion of a primitive Gromov--Witten
invariants and show that all Gromov--Witten invariants can be expressed in
terms of primitive invariants. In the case of the local invariants of a
nonsingular curve, this suggests a natural grading on the set of local
Gromov--Witten invariants. We will see that the \'etale invariants comprise
the first level of this grading.

\subsection{Inversion of the Gopakumar--Vafa formula}\label{subsec: inversion of GP}

Let $\beta \in H_{2} (X,\znums )$ be an indivisible class. Then the
Gopakumar--Vafa formula is:
\[
\sum _{g\geq 0}\sum _{d>0}N^{g}_{d\beta } (X)\lambda ^{2g-2}q^{d\beta } =
\sum _{g\geq 0}\sum _{d>0}n_{d\beta }^{g} (X)\sum
_{k>0}\tfrac{1}{k}\left(2\sin (\tfrac{k\lambda }{2})
\right)^{2g-2}q^{kd\beta}.
\]
Fix $n$ and look at the $q^{n\beta }$ terms on each side:
\[
\sum _{g\geq 0 }N_{n\beta }^{g} (X)\lambda ^{2g-2}=\sum _{g\geq 0}\sum _{d|n} n_{d\beta }^{g} (X)\tfrac{d}{n}\left(2\sin (\tfrac{n\lambda }{2d}) \right)^{2g-2}.
\]
Letting $s=n\lambda $ and multiplying the above equation by $n$ we find
\[
\sum _{g\geq 0}N_{n\beta }^{g} (X)n^{3-2g}s^{2g-2}=\sum _{d|n}\sum _{g\geq
0}n_{d\beta }^{g} (X)d\left(2\sin \tfrac{s}{2d} \right)^{2g-2}.
\]
Recall that M\"obius inversion says that if $f (n)=\sum _{d|n}g (d),$ then
$g (d)=\sum _{k|d}\mu (\frac{d}{k})f (k)$. Applying this to the above
equation (more precisely, to the coefficients of each term of the equation
separately), we obtain
\[
\sum _{g\geq 0}n_{d\beta }^{g} (X)d\left(2\sin \tfrac{s}{2d} \right)^{2g-2}
= \sum _{k|d}\mu (\tfrac{d}{k})\sum_{g\geq 0}N_{k\beta }^{g} (X)
(\tfrac{s}{k})^{2g-2}k.
\]
Letting $t=2\sin \frac{s}{2d}$ and dividing by $d$ we arrive at
\[
\sum _{g\geq 0}n_{d\beta }^{g} (X)t^{2g-2} = \sum _{g\geq 0}\sum _{k|d}\mu
(\tfrac{d}{k}) (\tfrac{d}{k})^{2g-3}N_{k\beta }^{g} (X)\left(2\arcsin
\tfrac{t}{2} \right)^{2g-2}.
\]
By interchanging $k$ and $d/k $ in the sum and restricting to the
$t^{2g-2}$ term of the formula we arrive at the following formula for the
BPS invariants.
\begin{proposition}\label{prop: general inversion of GV}
Let $\beta \in H_{2} (X,\znums )$ be an indivisible class, then the BPS
invariant $n_{d\beta }^{g} (X)$ is given by the following formula
\[
n_{d\beta }^{g} (X)=\sum _{g'=0}^{g}\sum _{k|d}\mu (k)k^{2g'-3} \alpha
_{g,g'} N^{g'}_{d\beta /k} (X)
\]     
where $\alpha _{g,g'}$ is the coefficient of $r^{g-g'}$ in the series 
\[
\left(\frac{\arcsin (\sqrt{r}/2)}{\sqrt{r}/2} \right)^{2g'-2}.
\]
In particular, $n_{d\beta }^{g} (X)$ depends on $N^{g'}_{d'\beta } (X)$ for
all $g'\leq g$ and all $d'$ dividing $ d$ such that $\mu (\frac{d}{d'})\neq
0$.
\end{proposition}
Note that the local BPS invariants are thus given by
\begin{equation}\label{eqn: inverted GV, formula for each coef}
n_{d}^{h} (g)=\sum _{h'=0} ^{h}\sum _{k|d}\mu (k)k^{2 (g+h')-3} \alpha
_{h+g,h'+g} N^{h'}_{d/k} (g),
\end{equation}
or in generating function form:
\begin{equation}\label{eqn: inverted GV, gen fnc form}
\sum _{d>0}\sum _{h\geq 0}n_{d}^{h} (g)t^{2 (g+h-1)}q^{d} =\!\!\! \sum _{k,n>0}\sum _{h\geq 0}\mu (k)N_{n}^{h} (g)k^{2 (g+h)-3}\left(2\arcsin \tfrac{t}{2} \right)^{2 (g+h-1)}q^{nk}.
\end{equation}

\subsection{Primitive Gromov--Witten invariants}\label{subsec: primitive GW
inv's} 

In this subsection, we formalize the observation that certain contributions
to the Gromov--Witten invariants of $X$ can be computed in terms of
Gromov--Witten invariants of the covering spaces of $X$. We use this to
reduce the computation of the \'etale invariants to the degree 1 invariants
(which have been previously computed by the second author \cite{Pa}).

\begin{definition}\label{defn: primitive GW invariants}
We say that a stable map $f\co C\to X$ is \emph{primitive} if $$f_{*}\co \pi _{1}
(C)\to \pi _{1} (X)$$ is surjective. Note that $\im (f_{*})\subset \pi _{1}
(X)$ is locally constant on the moduli space of stable maps. Let $\M _{g}
(X,\beta )_{G}$ be the component(s) consisting of maps $f$ with $\im
(f_{*})=G\subset \pi _{1} (X)$. In particular, $\M _{g} (X,\beta )_{\pi
_{1} (X)}$ consists of primitive stable maps. Define the \emph{primitive
Gromov--Witten invariants}, denoted $\hat{N}^{g}_{\beta } (X)$, to be the
invariants obtained by restricting $[\M _{g} (X,\beta )]^{vir}$ to the
primitive component $\M _{g} (X,\beta )_{\pi _{1} (X)}$.
\end{definition}

The usual Gromov--Witten invariants can be computed in terms of the
primitive invariants using the following observations. Let $\rho
\co \til{X}_{G}\to X$ be the covering space of $X$ corresponding to the
subgroup $G\subset \pi _{1} (X)$. Any stable map $$[f\co C\to X]\in \M _{g}
(X,\beta )_{G}$$ lifts to a (primitive) stable map $[\til{f}\co C\to
\til{X}_{G}]\in \M _{g} (\til{X}_{G},\til{\beta })_{G}$ for some
$\til{\beta }$ with $\rho _{*} (\til{\beta })=\beta $. Furthermore, this
lift is unique up to automorphisms of the cover $\rho \co \til{X}_{G}\to
X$. Conversely, any stable map in $\M _{g} (\til{X}_{G},\til{\beta })_{G}$
gives rise to a map in $\M _{g} (X,\beta )_{G}$ by composing with $\rho
$. Note that the automorphism group of the cover is $\pi _{1} (X)/N (G)$
where $N (G)$ is the normalizer of $G\subset \pi _{1} (X)$. If $G$ is
finite index in $\pi _{1} (X)$, then $\til{X}_{G}$ is compact and the
automorphism group of the cover is finite. This discussion leads to:
\begin{proposition}\label{prop: GW inv's of X in terms of primitive GW}
Fix $X$, $g$, and $\beta $. Suppose that for every stable map $[f:$\break
$C\to X]$
in $\M _{g} (X,\beta )$, the index $[\pi _{1} (X):f_{*} (\pi _{1} (C))]$ is
finite.  Then
\[
N_{\beta }^{g} (X)=\sum _{G} \sum _{\til{\beta }} \frac{1}{[\pi _{1} (X):N
(G)]}\hat{N}_{\til{\beta }}^{g} (\til{X}_{G})
\]
where the first sum is over $G\subset \pi _{1} (X)$ and the second sum is over
$\til{\beta }\in H_{2} (\til{X}_{G},\znums ) $ such that $ \rho _{*}
(\til{\beta })=\beta$.
\end{proposition}
\begin{remark}
In the case when $[\pi _{1} (X):G]=\infty $, $\til{X}_{G}$ will not be
compact and hence the usual Gromov--Witten invariants are not
well-defined. However, this technique sometimes can still be used to
compute the invariants (see \cite{BKL}). This technique originated in
\cite{Br-Le1} where it was used to compute multiple cover contributions of
certain nodal curves in surfaces.
\end{remark}

This technique is especially well-suited to the case of the local
invariants of a nonsingular genus $g$ curve. In this case, the image of the
fundamental group under a (non-constant) stable map always has finite
index. Furthermore, any degree $k$, complete, \'etale cover of a nonsingular
genus $g$ curve is a nonsingular curve of genus $k (g-1)+1$. Thus the formula in
Proposition~\ref{prop: GW inv's of X in terms of primitive GW} reduces to
\begin{equation}\label{eqn: local inv in terms of local primitive inv}
N_{n}^{h} (g)=\sum _{l|n}C_{l,g}\hat{N}_{n/l}^{h- (l-1) (g-1)} (l (g-1)+1)
\end{equation}
where $C_{k,g}$ is the number of degree $k$, connected, complete, \'etale
covers of a nonsingular genus $g$ curve, each counted by the reciprocal of the
number of automorphisms. In light of this formula, we can regard the
primitive local invariants $\hat{N}^{h}_{d} (g)$ as the fundamental
invariants. We encode these invariants into generating functions as
follows:
\[
\hat{F}_{k,g-1} (\lambda )=\sum _{h\geq 0}\hat{N}^{h}_{k} (g)\lambda ^{2
(g+h-1)}.
\]
Equation~\ref{eqn: local inv in terms of local primitive inv} can then be
written in generating function form as
\begin{align*}
\sum _{h\geq 0}\sum _{n>0}N^{h}_{n} (g)q^{n}t^{2 (g+h-1)} &= \sum _{h\geq 0}\sum _{k,l>0} C_{l,g}\hat{N}^{h- (l-1) (g-1)}_{k} (l (g-1)+1)q^{kl}t^{2 (g+h-1)}\\
&=\sum _{k,l>0}C_{l,g}\hat{F}_{k,l (g-1)}q^{kl}.
\end{align*}
We re-index and rearrange Equation~\ref{eqn: inverted GV, gen fnc form} below
\[
\sum _{d>0}\sum _{h\geq 0}n_{d}^{h} (g)t^{2 (g+h-1)}q^{d}=\!\!\sum
_{m>0}\tfrac{1}{m}\mu (m)\!\sum _{h\geq 0}\sum _{n>0}N_{n}^{h} (g) (q^{m})^{n}\left(2m\arcsin \tfrac{t}{2} \right)^{2 (g+h-1)}
\]
and then substitute the previous equation to arrive at the following
general equation for the local BPS invariants:
\[
\fbox{$\displaystyle \sum _{d>0}\sum _{h\geq 0}n_{d}^{h} (g)t^{2 
(g+h-1)}q^{d}=\!\!\!\!\sum_{m,k,l>0}\tfrac{1}{m}\mu  (m)C_{l,g}\hat{F}_{k,l (g-1)}
(2m\arcsin \tfrac{t}{2})^{2 (g+h-1)}q^{mkl}.
$}
\]
The unknown functions $\hat{F}_{k,l (g-1)}$ are graded by the two natural
numbers $k$ and $l$. The contribution in the above sum corresponding to
fixed $l $ and $k$ are from those stable maps that factor into a
composition of a degree $k$ primitive stable map and a degree $l$ \'etale
cover of $C_{g}$. Thus the \'etale BPS invariants (the type (i)
contributions) correspond exactly to restricting $k=1$ in the above
sum. Therefore we have
\[
\sum _{d>0}\sum _{h\geq 0}n_{d}^{h} (g)^{\et }t^{2
(g+h-1)}q^{d}=\sum_{m,l>0}\tfrac{1}{m}\mu (m)C_{l,g}\hat{F}_{1,l (g-1)}
(2m\arcsin \tfrac{t}{2})^{2 (g+h-1)}q^{ml}.
\]
Since a degree one map onto a nonsingular curve is surjective on the fundamental
group, it is primitive. The degree one local invariants were computed in
\cite{Pa}, the result can be expressed:
\begin{align*}
\hat{F}_{1,g-1}&=\sum _{h\geq 0}N_{1}^{h} (g)\lambda ^{2 (g+h-1)}\\
&=\left(4\sin ^{2}\tfrac{\lambda }{2} \right)^{(g-1)}
\end{align*}
and so
\[
\sum _{d>0}\sum _{h\geq 0}n_{d}^{h} (g)^{\et }t^{2
(g+h-1)}q^{d}=\sum_{m,l>0}\tfrac{1}{m}\mu (m)C_{l,g}\left(4\sin ^{2}
(m\arcsin \tfrac{t}{2}) \right)^{l (g-1)}q^{ml}.
\]
By the definition of $P_{m}$, we have 
\[
P_{m} (4\sin ^{2}\chi )=4\sin ^{2} (m\chi )
\]
and so letting $\chi =\arcsin (t/2)$ or equivalently $t=2\sin \chi $, we get
\[
\sum _{d>0}\sum _{h\geq 0}n_{d}^{h} (g)^{\et }t^{2
(g+h-1)}q^{d}=\sum_{m,l>0}\tfrac{1}{m}\mu (m)C_{l,g}\left(P_{m} (t^{2})) \right)^{l (g-1)}q^{ml}.
\]
Finally, by letting $y=t^{2}$ and re-indexing $m$ by
$k$, we get
\[
\sum _{d>0}\sum _{h\geq 0}n_{d}^{h} (g)^{\et }y^{h}q^{d}=\sum_{m,l>0}k\mu
(k)C_{l,g}P_{m} (y)^{l (g-1)}q^{ml},
\]
and so the formula in Theorem~\ref{thm: formula for etale BPS invariants}
is proved by comparing the $q^{d}$ terms.\qed

\section{Integrality of the \'etale BPS invariants}\label{sec: pf of
integrality}

In this section we show how the integrality of the \'etale BPS invariants
(Theorem \ref{thm: etale BPS invariants are integers}) follows from our formula
for them (Theorem \ref{thm: formula for etale BPS invariants}) and some
properties of the the polynomials $P_{l} (y)$ and the number of degree $k$
covers $C_{k,g}$.

The facts that we need concerning the polynomials $P_{l} (y)$ are the
following.

\begin{Plemmas}[Moll]\label{prop: formula for the poly P}
If $l\in \nnums $, then $P_{l} (y)$, defined by $P_{l} (4\sin ^{2}t)=4\sin
^{2} (lt)$, is given explicitly by
\[
P_{l} (y)=\sum _{a=1}^{l}-\frac{l}{a}{\binom{a-1+l}{2a-1}} (-y)^{a}.
\]
\end{Plemmas}

\begin{Plemmas}\label{lem: P is an integer poly}
If $l $ is a positive integer, then $P_{l} (y)$ is a polynomial with
integer coefficients.
\end{Plemmas}

\begin{Plemmas}\label{lem: Pab(y)=Pa(y)Pb(yPa(y))}
For any $\alpha $ and $\beta $ we have
\[
P_{\alpha \beta } (y)=P_{\beta } (P_{\alpha } (y)).
\]
\end{Plemmas}

\begin{Plemmas}\label{lem: mod p property of yPp(y)}
For $p$ a prime number and $b$ a positive integer, we have
\[
P_{p} (y)^{p^{l-1}b}\equiv y^{p^{l}b} \bmod p^{l}.
\]
\end{Plemmas}

We also will need some facts about $C_{k,g}$, the number of connected
\'etale covers.

\begin{Clemmas}\label{lem: a=A/k! is an integer}
Let $C $ be a nonsingular curve of genus $g$, let $S_{k}$ be the
symmetric group on $k$ letters, and define
\[
A_{k,g}=\# \Hom (\pi _{1} (C),S_{k}).
\]
Then 
\[
a_{k,g}=\frac{A_{k,g}}{k!}
\]
is an integer.
\end{Clemmas}

Note that $A_{k,g}$ is the number of degree $k$ (not necessarily connected)
\'etale covers of $C$ with a marking of one fiber. Thus $a_{k,g}$ is
the number of (not necessarily connected) \'etale covers each counted by
the reciprocal of the number of automorphisms. We remark that
Lemma~\ref{lem: a=A/k! is an integer} was essentially known to Burnsides.

\begin{Clemmas}\label{lem: gen fnc for C is the log of gen fnc of a}
Let $a_{k,g} $ be as above with $a_{0,g}=1$ by convention, then
\[
\sum _{k=1}^{\infty }C_{k,g}t^{k}=\log (\sum _{k=0}^{\infty }a_{k,g}t^{k}).
\]
\end{Clemmas}

\begin{Clemmas}\label{lem: c:=kC is an integer}
Define $c_{k,g}:=kC_{k,g}$. Then $c_{k,g}$ is an integer.
\end{Clemmas}

We remark that in general, $C_{k,g}$ is not an integer (see Table~3).

\begin{Clemmas}\label{lem: mod p properties of the c's}
Let $p$ be a prime number not dividing $k$ and let $l$ be a positive
integer. Then
\[
c_{p^{l}k,g}\equiv c_{p^{l-1}k,g} \bmod p^{l}.
\] 
\end{Clemmas}

We defer the proof of these lemmas to the subsections to follow and we
proceed as follows.

In light of Lemmas \ref{lem: P is an integer poly} and \ref{lem: c:=kC is
an integer}, we see from the formula in Theorem \ref{thm: formula for etale
BPS invariants} that $n_{d}^{h} (g)^{\et }\in \znums $ if and only if $\Xi
_{d,g}\equiv 0 \bmod d$, where
\[
\Xi _{d,g}=\sum _{k|d}\mu (k)c_{\frac{d}{k},g} P_{k} (y)^{\frac{d (g-1)}{k}}.
\]
Suppose that $p^{l}$ divides $d$ and that $p^{l+1}$ does not divide $d$ for
some prime number $p$.  For notational clarity, we will suppress the second
subscript of $c$ (which is always $g$) in the following calculation. Let
$a=d/p^{l}$; then we get
\begin{align*}
\Xi _{d,g}&=\sum _{k|a}\sum _{i=0}^{l}\mu (p^{i}k)\;c_{\frac{d}{p^{i}k}}\;
P_{p^{i}k} (y) ^{ \frac{d (g-1)}{p^{i}k}}\\
&=\sum _{k|a} \mu ({k})\; c_{\frac{p^{l}a}{k}}\;  P_{{k}} (y) ^{\frac{p^{l}a
(g-1)}{k}}\quad  -\mu ({k})\;c_{\frac{p^{l-1}a}{k}}\; P_{p{k}} (y)
^{\frac{p^{l-1}a (g-1)}{k}}.
\end{align*}
Let $\chi =P_{k} (y)$. Then by Lemma \ref{lem: Pab(y)=Pa(y)Pb(yPa(y))}
we have 
$P_{p{k}} (y)= P_{p} (\chi )$
and so 
\[
\Xi _{d,g}=\sum _{k|a}\mu  (k)\left\{c_{\frac{p^{l}a}{k}}\chi ^{\frac{p^{l}a (g-1)}{k}}-c_{\frac{p^{l-1}a}{k}}P_{p} (\chi )
^{\frac{p^{l-1}a (g-1)}{k}} \right\}.
\]
Then by Lemmas \ref{lem: mod p property of yPp(y)} and \ref{lem: mod p
properties of the c's} we have
\begin{align*}
\Xi _{d,g}&\equiv \sum _{k|a}\mu (k) 
\left\{c_{\frac{p^{l}a}{k}}\chi ^{\frac{p^{l}a (g-1)}{k}}-c_{\frac{p^{l}a}{k}}\chi ^{\frac{p^{l}a (g-1)}{k}} \right\} \bmod p^{l}\\
&\equiv 0 \bmod p^{l}
\end{align*}
and so $\Xi _{d,g}\equiv 0\bmod d$ and thus $n_{d}^{h} (g)^{\et }\in \znums
$. \qed

\subsection{Properties of the polynomials $P_{l} (y)$: the proofs of
Lemmas \ref{prop: formula for the poly P}--\ref{lem: mod p
property of yPp(y)}}\label{subsec: properties of the P polys}

This subsection is independent of the rest of the paper. We prove various
properties of the following family of power series:
\begin{definition}\label{defn: Pl(y)}
Let $\alpha \in \rnums $, we define the formal power series $P_{\alpha }
(y)$ by
\[
P_{\alpha } (y)=4\sin ^{2} (\alpha t)
\]
where 
\[
y=4\sin ^{2}t.
\]
Note that $P_{\alpha } (y)\in \rnums [[y]]$ since $\sin ^{2} (\alpha t)$ is
a power series in $t^{2}$ and $y (t)=4\sin
^{2}t=4t^{2}-\tfrac{4}{3!}t^{4}+\dots $ is an invertible power series in
$t^{2}$. (\textbf{Warning:} This definition differs from the one in
\cite{Br-RIMS} by a power of $y$.)
\end{definition}
\proof[Proof of Lemma \ref{lem: Pab(y)=Pa(y)Pb(yPa(y))}] 
This is immediate from the definition.
\qed\medskip 

\proof[Proof of Lemma \ref{prop: formula for the poly
P}] We prove the formula for $P_{l} (y)$ with $l\in \nnums $. This formula
and its proof was discovered by Victor Moll; we are grateful to him for
allowing us to use it.

From \cite{Whittaker-Watson} page 170 we can express $\sin ^{2} (lt)/\sin
^{2}t$ in terms of $\cos (2jt)$ for $1\leq j\leq l-1$ and from
\cite{Gradshteyn-Ryzhik} 1.332.3 we can in turn express $\cos (2jt)$ in
terms of $\sin ^{2}t$. Substituting, rearranging, and simplifying we arrive
a formula for the coefficients of $P_{l}$. Let $P_{l} (y)=\sum
_{n=1}^{l}-p_{n,l} (-y)^{n}$, then $p_{1,l}=l^{2}$ and for $l>1$,
\begin{equation}\label{eqn: first expression for coef of P}
p_{n,l}=\frac{1}{n-1}\sum _{j=n}^{l} (l-j+1) (j-1) {\binom{j+n-3}{j-n}}. 
\end{equation}
By standard recursion methods (see, for example, the book ``$A=B$''
\cite{A=B}) one can derive the identity for the binomial sum that transforms
the above expression for $p_{n,l}$ into the one asserted by the Lemma:
\begin{equation}\label{eqn: second formula for coef of P}
p_{n,l}=\frac{l}{n}{\binom{l+n-1}{2n-1}}.
\end{equation}
\qed\medskip

\proof[Proof of Lemma \ref{lem: P is an integer poly}] We need to show
that $p_{n,l}\in \znums $. By Equation \ref{eqn: second formula for coef of
P}, we have that $np_{n,l}\in \znums $ and by Equation \ref{eqn: first
expression for coef of P}, we have that $(n-1) p_{n,l}\in \znums $. Thus
$np_{n,l}- (n-1) p_{n,l}=p_{n,l}\in \znums$.\qed\medskip

Note that $-P_{l} (-y)$ has all positive integral coefficients.

\medskip
\proof[Proof of Lemma \ref{lem: mod p property of yPp(y)}]
To prove the lemma, clearly it suffices to prove that 
\[
P_{p} (y)^{p^{l-1}}\equiv y^{p^{l}}\bmod p^{l}
\]
for $p $ prime and $l\in N$.

For $n<p$, we have that  $p$ divides $p_{n,p}$ since 
\[
p_{n,p}=\frac{p}{n}\binom{p+n-1}{2n-1}
\]
and $n$ does not divide $p$ (except $n=1$). Noting that $p_{p,p}=1$ we have
\[
P_{p} (y)=y^{p}+pyf (y)
\]
for $f\in \znums [y]$. This proves the lemma for $l=1$. Proceeding by
induction on $l$, we assume the lemma for $l-1$ so that we can write
\[
P_{p} (y)^{p^{l-1}}=y^{p^{l-1}}+p^{l-1}g (y)
\]
where $g (y)\in \znums [y]$. But then
\begin{align*}
P_{p} (y) ^{p^{l}}&=\left(y^{p^{l-1}}+p^{l-1}g (y) \right)^{p}\\
&=y^{p^{l}}+\text{terms that $p^{l}$ divides}
\end{align*}
and so the lemma is proved. \qed

\subsection{Properties of the number of covers: the proofs of Lemmas
\ref{lem: a=A/k! is an integer}--\ref{lem: mod p properties
of the c's}}\label{subsec: properties of the C's}

In this subsection we prove the properties concerning the numbers $A_{k,g}$,
$a_{k,g}$, $C_{k,g}$, and $c_{k,g}$ that were asserted by the Lemmas.

We begin with a proposition from group theory due to M. Aschbacher:

\begin{proposition}[Aschbacher]\label{prop: Aschbacher's result}
Let $G$ be a finite group with conjugacy class\-es $C_{i}$, $1\leq i\leq
r$. Pick a representative $g_{i}\in C_{i}$; define
\[
b_{i,j,k}=|\{(g,h)\in C_{i}\times c_{j}:gh=g_{k} \}|
\]
and
\[
\beta _{i}=|\{(g,h)\in G\times G:[g,h]\in C_{i} \}|.
\]
Then
\[
\beta _{k}=|G|\cdot \sum _{i=1}^{r}b_{i,k,i}
\]
so, in particular, $|G|$ divides $\beta _{k}$. 
\end{proposition}

\proof
We use the notation $h^{-g}:=g^{-1}h^{-1}g$. For $(g,h)\in G\times G$,
\[
[g,h]=g^{-1}h^{-1}gh=h^{-g}h\in (h^{-1})^{G}\cdot h^{G}.
\]
Furthermore, $[x,h]=[y,h]$ if and only if $h^{-x}=h^{-y}$ if and only if
$xh^{-1}\in C_{G} (h)$, so
\[
\beta _{j,k}=|\{(g,h)\in G\times C_{j}:[g,h]=g_{k} \}|=|C_{G} (g_{j})|\cdot
b_{j',j,k}
\]
where $C_{j'}$ is the conjugacy class of $g_{j}^{-1}$. Of course
\[
\beta _{k}=|C_{k}|\sum _{j=1}^{r}\beta _{j,k}
\]
so 
\[
\beta _{k}=|C_{k}|\cdot \sum _{j=1}^{r}|C_{G} (g_{j})|b_{j',j,k}.
\]
Let
\[
\Omega_{j,k} =\{(g,h)\in C_{j'}\times C_{j}:gh\in C_{k} \}.
\]
Then
\begin{align*}
|\Omega _{j,k}|&=|C_{j}|\cdot |\{g\in C_{j'}:gg_{j}\in C_{k}
\}|\\
&=|C_{j}|\cdot |\{(g,h)\in C_{j'}\times C_{k}:g_{j}=g^{-1}h \}|\\
&=|C_{j}|b_{j,k,j}
\end{align*}
and similarly
\begin{align*}
|\Omega _{j,k}|&=|C_{k}|\cdot |\{(g,h)\in C_{j'}\times C_{j}:gh=g_{k}
\}|\\
&=|C_{k}|b_{j',j,k}
\end{align*}
so $|C_{k}|b_{j',j'k}=|C_{j}|b_{j,k,j}$. Therefore
\begin{align*}
\beta _{k}&=\sum _{j=1}^{r}|C_{G} (g_{j})|b_{j',j,k}|C_{k}|\\
&=\sum
_{j=1}^{r}|C_{G} (g_{j})|\cdot |C_{j}|b_{j,k,j}\\
&=|G|\sum
_{j=1}^{r}b_{j,k,j}
\end{align*}
which proves the proposition. \qed 
\medskip

\proof[Proof of Lemma \ref{lem: a=A/k! is an integer}] Recall that the
lemma asserts that $d!$ divides
\[
A_{k,g}=\# \Hom (\pi _{1} (C _{g}),S_{d}).
\]
For $x\in S_{d}$
let $c (x)$ denote the conjugacy class of $x$. We will prove, by induction
on $g$, that $d!$ divides the number of solutions $(x_{1},\dots
,x_{g},y_{1},\dots ,y_{g})$ to
\begin{equation}\label{eqn: [x1,y1]...[x2,y2]in c(z)}
\prod _{i=1}^{g}[x_{i},y_{i}]\in c (z)
\end{equation}
where $z$ is fixed. The lemma is then the special case where $z$ is the
identity.

The case of $g=1$ is Proposition \ref{prop: Aschbacher's result} where
$G=S_{d}$. 
For each fixed $r\in S_{d}$, the number of solutions to (\ref{eqn:
[x1,y1]...[x2,y2]in c(z)}) with 
\[
\prod _{i=1}^{g-1}[x_{i},y_{i}]=r
\]
is the
number of solutions to
\begin{equation}\label{eqn: w[xg,yg]inv=r}
w[x_{g},y_{g}]^{-1}=r
\end{equation}
as $w$ varies over $c (z)$ and $x_{g}$ and $y_{g}$ each vary over
$S_{d}$. The number of solutions to (\ref{eqn: w[xg,yg]inv=r}) depends only
on the conjugacy class of $r$ since if $q=srs^{-1}$, then (\ref{eqn:
w[xg,yg]inv=r}) holds if and only if
\[
sws^{-1}[sx_{g}s^{-1},sy_{g}s^{-1}]^{-1}=q
\]
holds. Thus the number of solutions to (\ref{eqn: [x1,y1]...[x2,y2]in
c(z)}) can be counted by summing up over $\{C_{i} \}$, the set of conjugacy
classes of $S_{d}$, the product of
\[
|\{(x_{g},y_{g},w)\in S_{d}\times S_{d}\times c (z):w[x_{g},y_{g}]^{-1}\in
C_{i} \}|
\]
with
\[
|\{(x_{1},\dots ,x_{g-1},y_{1},\dots y_{g-1})\in
(S_{d})^{2g-2}:\prod _{i=1}^{g-1}[x_{i},y_{i}]\in C_{i} \}|.
\]
By the induction hypothesis, this latter term is always divisible by $d!$,
thus the sum is also divisible by $d!$.\qed\medskip

\proof[Proof of Lemma \ref{lem: gen fnc for C is the log of gen fnc of
a}] 

$A_{k,g}$ is the number of $k$--fold (not necessarily connected), complete
\'etale covers of a nonsingular genus $g$ curve $C_{g}$ with a fixed labeling of
one fiber (the bijection is given by monodromy). Thus $a_{k,g}$ is the
number of such covers (without the label), each counted by the reciprocal
of the number of its automorphisms.

 The relationship between $a_{k,g}$, the total number of $k$--covers, and
$C_{k,g}$, the number of connected covers, is given by
\[
a_{k,g}=\sum _{\alpha = (1^{\alpha _{1}}2^{\alpha _{2}}\ldots)\in P
(k)}\frac{1}{\prod _{i\geq 1}\alpha _{i}!}\prod _{i\geq 1}C_{i,g}^{\alpha
_{i}}
\]
where $P (k)$ is the set of partitions of $k$ ($\alpha _{i}$ is the number
of $i$'s in the partition so $k=\sum i\alpha _{i}$). This formula is easily
obtained by considering how each cover breaks into a union of connected
covers (keeping track of the number of automorphisms).

Thus we have
\begin{align*}
\sum _{k=0}^{\infty }a_{k,g}t^{k}&=\sum _{\alpha _{1}\geq 0} \sum _{\alpha
_{2}\geq 0}\dots \prod _{i=1}^{\infty }\frac{1}{\alpha _{i}!}C_{i,g}^{\alpha
_{i}}t^{i\alpha
_{i}}\\
&=\sum _{\alpha \geq 0}\frac{1}{\alpha !}\left(\sum _{i=1}^{\infty
}C_{i,g}t^{i} \right)^{\alpha }\\
&=\exp (\sum _{i=1}^{\infty }C_{i,g}t^{i})
\end{align*}
and so 
\[
\sum _{k=1}^{\infty }C_{k,g}t^{k}=\log (\sum _{k=0}^{\infty }a_{k,g}t^{k})
\]
which proves the lemma.\qed\medskip 

\proof[Proof of Lemma \ref{lem: c:=kC is an integer}] Recall that the
lemma asserts that $c_{k,g}:=kC_{k,g}$ is an integer. From the previous lemma
we have:
\[
t\frac{d}{dt}\left(\sum _{k=1}^{\infty }C_{k,g}t^{k}
\right)=t\frac{d}{dt}\log\left(\sum _{k=0}^{\infty }a_{k,g}t^{k} \right)
\]
therefore
\[
\sum _{k=1}^{\infty }c_{k,g}t^{k}=\frac{\sum _{k=1}^{\infty
}ka_{k,g}t^{k}}{\sum _{k=0}^{\infty }a_{k,g}t^{k}}
\]
which implies
\[
la_{l,g}=\sum _{n=0}^{l-1}a_{n,g}c_{l-n,g}.
\]
Now $a_{0,g}=1$ and so we can obtain the $c$'s recursively from the $a$'s
and then induction immediately implies that $c_{l,g}\in \znums $. \qed
\medskip

\proof[Proof of Lemma \ref{lem: mod p properties of the c's}] We want to
prove that if $p $ is a prime number not dividing $k$ and $l$ a positive
integer, then
\[
c_{p^{l}k,g}\equiv c_{p^{l-1}k,g} \bmod p^{l}.
\]
We begin with two sublemmas: 
\begin{lem}\label{lem: mod p to l congruence for power of x+y}
Let $p $ be a prime, $l$ a positive number, and $x$ and $y$ variables, then
\[
(y+x)^{p^{l}}\equiv (y^{p}+x^{p})^{p^{l-1}} \bmod p^{l}.
\]
\end{lem}

\proof We use induction on $l$; the case $l=1$ is well known. By
induction, we may assume that there exists $\alpha \in \znums [x,y]$ such
that
\[
(y+x)^{p^{l-1}}= (y^{p}+x^{p})^{p^{l-2}}+\alpha p^{l-1}.
\]
Thus 
\begin{align*}
(y+x)^{p^{l}}&=\left((y^{p}+x^{p})^{p^{l-2}}+\alpha p^{l-1} \right)^{p}\\
&= (y^{p}+x^{p})^{p^{l-1}}+\text{terms divisible by $p^{l}$}
\end{align*}
which proves the sublemma.

\begin{lem}\label{lem: p to the l-k divides p to the l-1 choose k}
Let $1\leq k\leq l-1$ and let $p$ be prime. Then
$p^{l-k}$ divides ${\binom{p^{l-1}}{k}}$.
\end{lem}

\proof Recall Legendre's formula for $v_{p} (m!)$, the number of
$p$'s in the prime decomposition of $m!$:
\[
v_{p} (m!)=\frac{m-S_{p} (m)}{p-1}
\]
where $S_{p} (m)$ is the sum of the digits in the base $p$ expansion of $m$. 

Let $(a_{l-2},\dots ,a_{0})$ and $(b_{l-2},\dots ,b_{0})$ be base $p$
expansions of $k$ and $p^{l-1}-k$ respectively, then a simple calculation
yields:
\[
v_{p} ({\binom{p^{l-1}}{k}})=\frac{1}{p-1} (\sum a_{i}+\sum b_{i}-1).
\]
Let $n=v_{p} (k)$ so that $a_{n}$ is the first non-zero digit of $k=
(a_{l-2},\dots ,a_{0}) $. Now addition in base $p$ gives $(a_{l-2},\dots
a_{0})+ (b_{l-2},\dots ,b_{0})= (1,0,\dots ,0)$ so we have that
$b_{0}=b_{1}=\dots =b_{n-1}=0$, $b_{n}=p-a_{n}$, and $b_{i}=p-1-a_{i}$ for
$n+1\leq i\leq l-2$. Thus we see that 
\[
\sum a_{i}+\sum b_{i}-1= (l-1-n) (p-1)
\]
and so, observing that $k\geq n+1$, we have
\[
v_{p} ({\binom{p^{l-1}}{k}})=l-1-n\geq l-k
\]
which proves the sublemma.

Now let $a (t)=\sum _{i=1}^{\infty }a_{i,g}t^{i}$ so that Lemma \ref{lem: gen
fnc for C is the log of gen fnc of a} can be written
\[
\sum _{i=1}^{\infty }c_{i,g}\frac{t^{i}}{i}=\log (1+a (t)).
\]
Thus we have
\begin{align*}
c_{p^{l}k,g}&=p^{l}k\Coeff{t^{p^{l}k}} \{\log (1+a (t)) \}\\
c_{p^{l-1}k,g}&=p^{l-1}k\Coeff{t^{p^{l}k}} \{\log (1+a (t^{p})) \}
\end{align*}
and so
\begin{align*}
c_{p^{l}k,g}-c_{p^{l-1}k,g} &=k\Coeff{t^{p^{l}k}} 
\left\{\log\left(\frac{(1+a (t))^{p^{l}}}{(1+a (t^{p}))^{p^{l-1}}} \right) \right\}\\
&=k\Coeff{t^{p^{l}k}}\left\{Q (t)+\frac{Q^{2} (t)}{2}+\frac{Q^{3}
(t)}{3}+\dots \right\}
\end{align*}
where $Q \in t\znums [[t]]$ is defined by 
\[
\frac{(1+a (t))^{p}}{1+a (t^{p})}=1-Q (t).
\]
To prove Lemma \ref{lem: mod p properties of the c's} it suffices to prove
that $Q (t)\equiv 0 \bmod p^{l}$ since then $Q+Q^{2}/2+Q^{3}/3+\dots \in
\znums _{(p)}[[t]]$ and $Q+Q^{2}/2+Q^{3}/3+\dots\equiv 0 \bmod p^{l} $
which then proves that $c_{p^{l}k,g}-c_{p^{l-1}k,g}\equiv 0 \bmod p^{l}$.

Thus we just need to show that 
\[
(1+a (t))^{p^{l}}\equiv (1+a (t^{p}))^{p^{l-1}}\bmod p^{l}.
\]
From Lemma \ref{lem: mod p to l congruence for power of x+y} we have
\begin{align*}
(1+a (t))^{p^{l}}&\equiv (1+a (t)^{p})^{p^{l-1}}\bmod p^{l}\\
&\equiv (1+a (t^{p})+pf (t))^{p^{l-1}}\bmod p^{l}\\
&\equiv (1+a (t^{p}))^{p^{l-1}}+\sum _{k=1}^{p^{l-1}}
{\binom{p^{l-1}}
{k}}p^{k}f (t)^{k} (1+a (t^{p}))^{p^{l-1}-k} \bmod p^{l}.
\end{align*}
By Lemma \ref{lem: p to the l-k divides p to the l-1 choose k}, $p^{l}$
divides all the terms in the sum and thus Lemma \ref{lem: mod p properties
of the c's} is proved. \qed

\section{The Grothendieck--Riemann--Roch calculation}\label{sec: the GRR computation}

In order to prove Theorem~\ref{thm: GRR computation}, we will apply the
Grothendieck--Riemann--Roch formula to the morphism $\pi\co U \rightarrow
\til{M}$ of nonsingular stacks. Here $U$ is the universal curve; see
Subsection~\ref{subsec: type ii contributions} for the definition of
$\til{M}$. The first step is to compute the relative Todd class of the
morphism $\pi$---that is:
$$Td(\pi)= Td(T_U)/Td(\pi^*(T_{\til{M}})).$$ As the singularities of the
morphism $\pi$ occur exactly at the nodes of the universal curve (and the
deformations of the $1$--nodal map surject onto the versal deformation space
of the node), we may use a formula derived by D. Mumford for the relative
Todd class \cite{Mumford} (c.f. \cite{Fa-Pa} Section~1.1).

Let $S\subset U$ denote the (nonsingular) substack of nodes.
$S$ is of pure codimension 2. 
There is canonical double cover of $S$,
$$\iota\co Z \rightarrow S$$
obtained by ordering the branches of the node.
$Z$ carries two natural line bundles: the cotangent lines on the
first and second branches. Let $\psi_+$, $\psi_-$ denote
the Chern classes of these line bundles in $H^2(Z,\mathbb{Q})$.
Let $K=c_1(\omega_\pi)\in H^2(U,\mathbb{Q})$.
Mumford's formula is:
$$Td(\pi)= \frac{K}{e^K-1}+ \frac{1}{2}\iota_* \Big(
\sum_{l=1}^{\infty}\frac{B_{2l}}{(2l)!} \frac{\psi_+^{2l-1}+\psi_-^{2l-1}}
{\psi_+ +\psi_-} \Big).$$
Since $U$ is a threefold and $S$ is a curve, we find:
\begin{equation}
\label{qqq}
Td(\pi)= 1-\frac{K}{2}+\frac{K^2}{12} + \frac{[S]}{12}.
\end{equation}
Let $\gamma_0+ \gamma_1+\gamma_2 \in H^*(\til{M},\mathbb{Q})$ denote the
cohomological $\pi$ push-forward of $Td(\pi)$.
The evaluations:
\begin{equation}
\label{gamun}
\gamma_0= -g-(d-1)(g-1), \ \ \gamma_1= \pi_*\Big(
\frac{K^2+ [S]}{12} \Big),
\ \ \gamma_2=0,
\end{equation}
follow from equation (\ref{qqq}).

The Grothendieck--Riemann--Roch formula determines
the Chern character of the $\pi$ push-forward:
$$ch(R^\bullet \pi_* f^*(\mathcal{O}_C)) = \pi_*(ch(f^* (\mathcal{O}_C)) \cdot
Td(\pi)).$$ 
The right side is just $\gamma_0+\gamma_1+\gamma_2$.
By GRR again,
\begin{equation}
\label{rrr}
ch(R^\bullet \pi_* f^*(\omega_C)) = \pi_*(ch(f^*(\omega_C))\cdot Td(\pi)).
\end{equation}
We may express the right side as  
$$\til{\gamma}_0+\til{\gamma}_1+
\til{\gamma_2} \in H^*(\til{M},\mathbb{Q})$$
by the following formulas:
\begin{equation}
\label{gamti}
\til{\gamma}_0= g-2-(d-1)(g-1), \ \ 
\til{\gamma}_1= \pi_*\Big(
\frac{K^2+ [S]}{12} - \frac{ K \cdot W}{2}\Big),
\end{equation}
$$\til{\gamma}_2 = \pi_*\Big(\frac{K^2\cdot W +[S]\cdot W}{12}\Big).$$
These equations are obtained by simply expanding (\ref{rrr}) where we use
the notation:
\[
W:=f^*(c_1(\omega_C)).
\]
The Chern characters of $R^\bullet\pi_* f^*(\mathcal{O}_C\oplus
\omega_C)$ determine 
the classes of the expression:
$$c( R^\bullet\pi_* f^*(\mathcal{O}_C\oplus
\omega_C)[1]).$$
A direct calculation shows:
\begin{equation}
\label{fred}
\int_{\til{M}} c( R^\bullet\pi_* f^*(\mathcal{O}_C\oplus
\omega_C)[1]) = \int_{\til{M}}
\til{\gamma}_2 + \frac{\gamma_1^2+\til{\gamma}_1^2}{2}
+\gamma_1\til{\gamma_1}.
\end{equation}
Therefore,
our next step is to compute 
the intersections of the $\gamma$ and $\til{\gamma}$
classes
in $H^4(\til{M},\mathbb{Q})$.

\subsection{$\Sym^2(C)$} Let $\Sym^2(C)$ be the symmetric product of
$C$. $\Sym^2(C)$ is a nonsingular scheme.  There is a canonical branch
morphism 
\[
\mu\co  \til{M} \rightarrow \Sym^2(C)
\]
which associates the branch
divisor to each point $[f\co D \rightarrow C] \in \til{M}$ (see
\cite{Fantechi-Pandharipande}).  The degree of the morphism $\mu$ is
$C_{d,g}$.  We will relate the required intersections in $\til{M}$ to the
simpler intersection theory of $\Sym^2(C)$.

Let $L \in H^2(\Sym^2(C),{\mathbb{Q}})$ 
denote the divisor class corresponding to the
subvariety: 
$$L_p= \{ (p,q) \ | \ q \in C \ \}.$$
Let $\Delta$ denote the diagonal divisor class of $\Sym^2(C)$.
It is easy to compute the products:
$$L^2=1, \ \  L\cdot \Delta = 2, \ \ \Delta^2 = 4-4g$$
in $\Sym^2(C)$.

\subsection{$R$, $S$, and $T$}
An analysis of the ramification of the universal map 
$f\co U \rightarrow C$ is required
to relate the integrals (\ref{fred}) over $\til{M}$ 
to the intersection theory of $\Sym^2(C)$.
Consider first the maps:
$$U \stackrel{\alpha}{\rightarrow} \til{M} \times C
\stackrel{\beta}{\rightarrow} \til{M}$$
where
$\alpha=(\pi,f)$ and $\beta$ is the projection onto
the first factor.
Let 
$$R\subset U,$$
$$B \subset \til{M}\times C,$$ 
denote universal ramification and branch loci respectively.
Certainly,
\begin{equation}
\label{qwe}
\alpha_*([R])= [B]
\end{equation}
as the $\alpha$ restricts to a birational morphism from $R$ to $B$.
By the Riemann--Hurwitz correspondence, we find:
\begin{equation}
\label{rrttrr}
K= W + [R] \in H^2(U,\mathbb{Q}).
\end{equation}
After taking the square of this equation
 and pushing forward via $\alpha$, we find the equation
\begin{equation}
\label{ned}
\alpha_*(K^2)= 2[B]\cdot c_1(\omega_C)+ \alpha_*([R]^2)
\end{equation}
holds on $\til{M}\times C$.

The term $\alpha_*([R]^2)$ in (\ref{ned}) may be determined by the
following considerations.  The line bundle $\omega^*_\pi|_R$ is naturally
isomorphic to ${\mathcal{O}}_U(R)|_R$ at each point of $R$ not contained in
the locus of nodes $S$ or the locus of double ramification points $T$.  We
will use local calculations to show that the coefficients of $[S]$ and
$[T]$ are 1 in the following equation:
\begin{equation}
\label{ddffgg}
[R]^2= - K\cdot [R] + [S] + [T].
\end{equation}
To compute the coefficient of $[S]$ it suffices to study the local family
$\pi \co U_{loc}\to \cnums $ given by $U_{loc}=\{(x,y,t)\in \cnums ^{3}: xy=t
\}$ with the maps $f (x,y,t) =x+y$ and $\pi (x,y,t)=t$. For the coefficient
of $[T]$, we note $\omega_\pi^*$ is the $\pi$--vertical tangent bundle of
$U$ on $T$. Near $T$, $R$ is a double cover of $\til{M}$ with simple
ramification at $T$. Hence, the natural map on $R$ near $T$:
$$\omega_\pi^*|_R \rightarrow  {\mathcal{O}}_U(R)|_R$$
has a zero of order 1 along T. The coefficient of 
$[T]$ in (\ref{ddffgg}) is thus 1.
We may rewrite (\ref{ddffgg}) using (\ref{rrttrr})
$$ [R]^2 = \frac{-W \cdot [R] + [S] +[T]}{2},$$
which will be substituted in (\ref{ned}).

The final equation for $\alpha_*(K^2)$ using the above
results is:
\begin{equation}
\label{asd} 
\alpha_*(K^2)= \frac{3}{2} 
[B]\cdot c_1(\omega_C)+ \alpha_*\Big(\frac{[S]+[T]}{2}\Big).
\end{equation}
Note
the branch divisor $B$ is simply the $\mu$ pull-back of the
universal family 
$$B_S \subset \Sym^2(C) \times C.$$
Let $\beta_S$ denote the projection of  $\Sym^2(C) \times C$
to the first factor.
Applying $\beta_*$ to (\ref{asd})
and using the $\mu$ pull-back structure of $B$, we find:
$$\pi_*(K^2)=\beta_*\alpha_*(K^2)= \mu^*\beta_{S*}\big(
\frac{3}{2} [B_S]\cdot c_1(\omega_C)) + \pi_*\Big(\frac{[S]+[T]}{2}\Big).$$
A simple calculation in $\Sym^2(C)\times C$ then
yields:
$$\beta_{S*}( [B_S]\cdot c_1(\omega_C)) = (2g-2)L,$$
We finally arrive at the central equation: 
\begin{equation}
\label{newt}
\pi_*(K^2)= \mu^*\big(\frac{3}{2}(2g-2)L\big) + \pi_*\Big(
\frac{[S]+[T]}{2}\Big).
\end{equation}
Equation (\ref{newt}) will be used to 
transfer intersections on $\til{M}$ to $\Sym^2(C)$.

\subsection{Proof of Theorem~\ref{thm: GRR computation}}
We will calculate all terms on the right side of 
integral equation: 
\begin{equation}
\label{fredddd}
\int_{\til{M}} c( R^\bullet\pi_* f^*(\mathcal{O}_C\oplus
\omega_C)[1]) = \int_{\til{M}}
\til{\gamma}_2 + \frac{\gamma_1^2+\til{\gamma}_1^2}{2}
+\gamma_1\til{\gamma_1}.
\end{equation}
Consider first the class $\til{\gamma}_2$. By equation (\ref{gamti}),
\begin{equation}
\label{wer}
\int_{\til{M}} \til{\gamma}_2 =
\int_{\til{M}} \frac{ \pi_*(K^2\cdot W) + [S]\cdot W}{12}.
\end{equation}
The first summand on the right may be computed from the
relation:
$$\pi_*(K^2\cdot W) = \frac{[S]\cdot W + [T]\cdot W}{2}.$$
The definition of the Hurwitz numbers $D^*_{d,g}$ and
$D^{**}_{d,g}$ imply:
$$[S]\cdot W = D^*_{d,g}(2g-2),$$
$$[T]\cdot W = D^{**}_{d,g}(2g-2).$$
Using the above formulas, we find:
$$\int_{\til{M}} \til{\gamma}_2 = \frac{g-1}{12} \Big(
 3 D^*_{d,g}+ D^{**}_{d,g} \Big).$$
For the quadratic terms involving $\gamma_1$ and $\til{\gamma}_1$
in equation (\ref{fredddd}), we will need to compute 
several integrals.
The first two integrals are:
$$
\int_{\til{M}} \pi_*([S])^2 = (4-4g) D^*_{g,d}, \ \
\int_{\til{M}} \pi_*([T])^2 = \frac{4-4g}{3} D^{**}_{g,d}.
$$
Both equations require a study of the local geometry of the morphism
$\mu$. As $\mu$ is \'etale at the points of $\pi(S)$, the self-intersection
of the curve $\pi(S)$ is simply $\Delta^2\cdot D^*_{d,g}$.  As $\mu$ has
double ramification at the points of $\pi(T)$, the self-intersection of the
curve $\pi(T)$ is one third of $\Delta^2 \cdot D^{**}_{d,g}$ (see
\cite{Harris-Morrison}).  The integral :
$$\int_{\til{M}} \pi_*(K^2)^2 = 
\frac{9}{4}(2g-2)^2 D_{d,g}
+(5g-5)D^*_{d,g}+ \frac{17g-17}{3}D^{**}_{d,g}.$$
then 
follows easily from (\ref{newt}).

Next, the integral
$$\int_{\til{M}} \pi_*(K^2) \cdot \pi_*([S]) =
(4g-4)D^*_{g,d}$$
follows from the intersection theory of $\Sym^2(C)$ and the
definition of the Hurwitz numbers.

Finally, as $\pi_*(K\cdot W) = \mu^*((2g-2) L)$,
the remaining integrals:
$$
\int_{\til{M}} \pi_*(K^2) \cdot \pi_*(K\cdot W) =
\frac{3}{2}(2g-2)^2 D_{d,g} + (2g-2)(D^*_{d,g}+D^{**}_{d,g}),$$
$$ \int_{\til{M}} \pi_*([S]) \cdot \pi_*(K\cdot W) =
(4g-4)D^*_{g,d}, $$
$$\int_{\til{M}}\pi_*(K\cdot W) ^2 = (2g-2)^2 D_{d,g},$$
are easily obtained.

The final formula for Theorem~\ref{thm: GRR computation} is now obtained
from the above integral equations together with (\ref{gamun}),
(\ref{gamti}), and (\ref{fredddd}):
$$
\int _{\til{M}}c (R^{\bullet }\pi _{*}f^{*}
(\mathcal{O}_{C}\oplus \omega _{C})[1])
 =
\frac{g-1}{8} \Big( (g-1)D_{d,g} - D_{d,g}^* -\frac{1}{27}
D_{d,g}^{**}\Big).$$

\appendix \section{Appendix: Tables of numbers}\label{appendix}

\small

The tables in this appendix list the values of the invariants
studied in the paper in the first few cases:
\begin{enumerate}
\item [(1)] The \'etale BPS invariants $n_{d}^{h} (g)^{\et }$ (for small values
of $d$, $g$, and $h$) as given by Theorem~\ref{thm: formula for etale BPS
invariants}.
\item [(2)] The full local BPS invariants $n_{d}^{h} (g)$ (again for small
values of $d$, $g$, and $h$) in the range where they are known as given by
Lemma~\ref{lem: local BPS invs in known range (including GRR)} --- question
marks where they are unknown.
\item [(3)] The various Hurwitz numbers that arise. 
\end{enumerate}
The Hurwitz numbers were computed from first principles and recursion when
possible (see \cite{Br-RIMS} for example), and by a naive computer
program elsewhere. The Hurwitz numbers that were beyond either of these
methods are left as variables in the tables. Note that by Lemma~\ref{lem:
gen fnc for C is the log of gen fnc of a}, the rational numbers $C_{d,g}$
can be expressed in terms of the integers $a_{d,g}$; it is easy to write a
computer program that computes the $a_{d,g}$'s (albeit slowly).

If the \'etale BPS invariants do indeed arise from corresponding
component(s) in the D--brane moduli space (see Remark~\ref{rem: D-bran
moduli space should have ``etale'' component}) then the horizontal rows of
the tables for the \'etale invariants should be the coefficients of the
$\mathfrak{sl_{2}\oplus sl_{2}}$ decomposition of the cohomology of that
space. So for example, the zeros in the beginning of the $n_{4}^{h}
(g)^{\et }$ table suggest that this space factors off (cohomologically)
a torus of dimension $2g-1$ (as oppose to the torus factor of dimension $g$
for the other cases).

\begin{table}[ht!]\tiny
\begin{align*}
&
\begin{array}{||c||c|c|c|c|c|c|c|c|c|c||}
\hline \hline 
n_{2}^{h} (g)^{\et }&h=0  &h=1&h=2&h=3&h=4&h=5&h=6&h=7&h=8&h=9 \\ \hline 
g=2 &-2 &8 &0 &0 &0 &0 &0 &0 &0 &0 \\ \hline
g=3 &-8 &4 &31 &0 &0  &0 &0 &0 &0&0 \\ \hline  
g=4 &-32 &24 &-6 &128 &0  &0 &0 &0 &0&0 \\ \hline
g=5 &-128 &128 &-48 &8 &511 &0 &0 &0 &0 &0 \\ \hline
g=6 &-512 &640 &-320 &80 &-10 &2048 &0 &0 &0 &0 \\ \hline
g=7 &-2048 &3072 &-1920 &640 &-120 &12 &8191 &0 &0 &0 \\ \hline  
\end{array}
\\
&\\
&
\begin{array}{||c||c|c|c|c|c|c|c|c||}
\hline \hline n_{3}^{h} (g)^{\et } &h=0&h=1&h=2&h=3&h=4&h=5&h=6&h=7 \\ \hline
 g=2 &-3 &2 &73 &0 &0 &0 &0 &0 \\ \hline 
g=3 &-27 &36 &-18 &4 &2641 &0 &0 &0 \\ \hline 
g=4 &-243 &486 &-405 &180 &-45 &6 &93913 &0 \\ \hline
g=5 &-2187 &5832 &-6804 &4536 &-1890 &504 &-84 &8 \\ \hline 
g=6 &-19683 &65610 &-98415 &87480 &-51030 &20412 &-5670 &1080 \\ \hline 
g=7 &-177147 &708588 &-1299078 &1443420 
&-1082565 &577368 &-224532 &64152 \\ \hline 
\end{array}
\\
&\\
&
\begin{array}{||c||c|c|c|c|c|c|c|c||}
\hline \hline n_{4}^{h} (g)^{\et } &h=0&h=1&h=2&h=3&h=4&h=5&h=6&h=7 \\ \hline
 g=2  &0 &-60 &30 &1315 &0 &0 &0 &0 \\ \hline 
g=3  &0 &0 &-4032 &4032 &-1512 &252 &689311 &0 \\ \hline 
g=4 &0&0&0&-261120 &391680 &-244800 &81600 &-15300  \\ \hline
g=5 &0&0&0&0&-16760832&33521664&-29331456&14665728 \\ \hline
\end{array}
\\
&\\
&
\begin{array}{||c||c|c|c|c|c|c|c|c||}
\hline \hline n_{5}^{h} (g)^{\et } &h=0&h=1&h=2&h=3&h=4&h=5&h=6&h=7 \\ \hline
 g=2 &-5 &10 &-7 &2 &-1935+a_{5,2} &0 &0 &0 \\ \hline 
g=3 &-125 &500 &-850 &800 &-455 &160 &-34 &4 \\ \hline 
g=4 &-3125 &18750 &-50625 &81250 &-86250 &63750 &-33625 &12750 \\ \hline
g=5 &-78125 &625000 &-2312500 &5250000 &-8181250 &9275000 &-7910000 &5175000 \\ \hline 
\end{array}
\end{align*}
\caption{The \'etale BPS invariants $n_{d}^{h} (g)^{\et }$ for small $d$,
$g$, and $h$.}
\end{table}

\begin{table}[ht!]\tiny
\begin{align*}
&
\begin{array}{||c||c|c|c|c|c|c|c|c|c|c||}
\hline \hline 
n_{2}^{h} (g)&h=0  &h=1&h=2&h=3&h=4&h=5&h=6&h=7&h=8&h=9 \\ \hline 
g=2 &-2 &8 &0 &? &? &? &? &? &? &? \\ \hline
g=3 &-8 &4 &31 &8 &?  &? &? &? &?&? \\ \hline  
g=4 &-32 &24 &-6 &128 &96  &? &? &? &?&? \\ \hline
g=5 &-128 &128 &-48 &8 &511 &768 &? &? &? &? \\ \hline
g=6 &-512 &640 &-320 &80 &-10 &2048 &5120 &? &? &? \\ \hline
g=7 &-2048 &3072 &-1920 &640 &-120 &12 &8191 &30720 &? &? \\ \hline  
\end{array}
\\
&\\
&
\begin{array}{||c||c|c|c|c|c|c|c|c||}
\hline \hline n_{3}^{h} (g) &h=0&h=1&h=2&h=3&h=4&h=5&h=6&h=7 \\ \hline
g=2 &-3 &2 &73 &50 &? &? &? &? \\ \hline 
g=3 &-27 &36 &-18 &4 &2641 &9604 &? &? \\ \hline 
g=4 &-243 &486 &-405 &180 &-45 &6 &692352 &836310 \\ \hline
g=5 &-2187 &5832 &-6804 &4536 &-1890 &504 &-84 &8 \\ \hline 
g=6 &-19683 &65610 &-98415 &87480 &-51030 &20412 &-5670 &1080 \\ \hline 
g=7 &-177147 &708588 &-1299078 &1443420 
&-1082565 &577368 &-224532 &64152 \\ \hline 
\end{array}
\\
&\\
&
\begin{array}{||c||c|c|c|c|c|c|c|c||}
\hline \hline n_{4}^{h} (g) &h=0&h=1&h=2&h=3&h=4&h=5&h=6&h=7 \\ \hline
 g=2  &0 &-60 &30 &? &? &? &? &? \\ \hline 
g=3  &0 &0 &-4032 &3520 &? &? &? &? \\ \hline 
g=4 &0&0&0&-261120 &367104 &? &? &?  \\ \hline
g=5 &0&0&0&0&-16760832&32735232&?&? \\ \hline
\end{array}
\\
&\\
&
\begin{array}{||c||c|c|c|c|c|c|c|c||}
\hline \hline n_{5}^{h} (g) &h=0&h=1&h=2&h=3&h=4&h=5&h=6&h=7 \\ \hline 
g=2 &-5 &10 &-7 &2 &-1935+a_{5,2} &* &? &? \\ \hline 
g=3 &-125 &500 &-850 &800 &-455 &160 &-34 &4 \\ \hline 
g=4 &-3125 &18750 &-50625 &81250 &-86250 &63750 &-33625 &12750 \\ \hline 
g=5 &-78125 &625000 &-2312500 &5250000 &-8181250 
&9275000 &-7910000 &5175000 \\ \hline
\end{array}
\end{align*}
\caption{The local BPS invariants $n_{d}^{h} (g)$ for small $d$,
$g$, and $h$. Note: the value of $*$ in the above table is $\tfrac{1}{8}
(D_{5,2}-D_{5,2}^{*}-\tfrac{1}{27}D^{**}_{5,2})$.}
\end{table}

\begin{table}[ht!]\tiny
\begin{align*}
&
\begin{array}{||c||c|c|c|c|c|c||}
\hline \hline 
C_{d,g}  &g=1&g=2&g=3&g=4&g=5&g=6 \\ \hline 
d=2 &3/2 &15/2 &63/2 &255/2  &1023/2  &4095/2    \\ \hline
d=3 &4/3 & 220/3&7924/3 &281740/3  &10095844/3  &362968060/3  \\ \hline
d=4 &7/4 &5275/4 &2757307/4 &a_{4, 4} - 408421/4  &a_{4, 5} - 13985413/4  &a_{4, 6} - 492346021/4  \\ \hline
\end{array}
\\
&\\
&
\begin{array}{||c||c|c|c|c|c|c|c||}
\hline \hline 
D_{d,g}  &g=1&g=2&g=3&g=4&g=5&g=6&g=7 \\ \hline 
d=2 & 2&8 &32 &128  & 512 &2048  &8192  \\ \hline
d=3 &16 &640 &23296 &839680  &30232576 &1088389120  &39182073856  \\ \hline
\end{array}
\\
&\\
&
\begin{array}{||c||c|c|c|c|c|c|c||}
\hline \hline 
D^{*}_{d,g}  &g=1&g=2&g=3&g=4&g=5&g=6&g=7 \\ \hline 
d=2 & 2&8 &32 &128  & 512 &2048  &8192   \\ \hline
d=3 &7 &235 &7987 &281995  &10096867  &362972155  &13062280147   \\ \hline
\end{array}
\\
&\\
&
\begin{array}{||c||c|c|c|c|c|c|c||}
\hline  \hline 
D^{**}_{d,g}  &g=1&g=2&g=3&g=4&g=5&g=6&g=7 \\ \hline 
d=2 & 0&0 &0 &0  &0  &0  &0   \\ \hline
d=3 & 3&135 &5103 &185895  &6711903  &241805655  &8706597903   \\ \hline
\end{array}
\end{align*}
\caption{The Hurwitz numbers $C_{d,g}$, $D_{d,g}$, $D^{*}_{d,g}$, and
$D_{d,g}^{**}$ for small $d$ and $g$.}
\end{table}

\end{document}